\documentclass[10pt]{article}
\usepackage[utf8]{inputenc}
\usepackage[english]{babel}
\usepackage{amsmath,amsthm,amsfonts}
\usepackage{graphicx}
\usepackage{bm}
\usepackage{tikz}
\usepackage{multirow}

\newcommand{\grad}{\mathop{\rm grad}\nolimits}
\renewcommand{\div}{\mathop{\rm div}\nolimits}

\newcommand{\norm}[1]{\left\|#1\right\|}

\newcommand{\seminorm}[1]{\left|#1\right|}

\begin{document}
\title{
A Generalized Multiscale Finite Element Method for Poroelasticity Problems I: Linear Problems
\thanks{This work was supported by RFBR (project N 15-31-20856)}}

\author{
Donald L. Brown
\thanks{
Institute for Numerical Simulation University of Bonn, 
Wegeler Strasse 6, 53115 Bonn, Germany
Email: {\tt donaldbrowdr@gmail.com}.
}
\and
Maria Vasilyeva
\thanks{ Institute of Mathematics and Informatics, 
North-Eastern Federal University, 
Yakutsk, Republic of Sakha (Yakutia), Russia, 677980 \& 
Institute for Scientific Computation, 
Texas A\&M University, College Station, 
TX 77843Email: {\tt vasilyevadotmdotv@gmail.com}.
}
}

\maketitle

\begin{abstract}
In this paper, we consider the numerical solution of poroelasticity problems that are of Biot type and develop a general algorithm for solving coupled systems. 
We discuss the challenges associated with mechanics and flow problems in heterogeneous media. 
The two primary issues being the multiscale nature of the media and the solutions of the  fluid and mechanics variables  traditionally developed with  separate grids and methods. 
For the numerical solution we develop and implement a Generalized Multiscale Finite Element Method (GMsFEM) that solves problem on a coarse grid by constructing local multiscale basis functions. 
The procedure begins with  construction of multiscale bases for both  displacement and pressure in each coarse block.  
Using a snapshot space and local spectral problems, we construct a basis of reduced dimension. Finally, after multiplying by a multiscale partitions of unity, the multiscale basis is constructed in the offline phase and 
the coarse grid problem then can be solved for arbitrary forcing and boundary conditions. 
We implement this algorithm on two heterogenous media and compute error  between the multiscale solution  with the fine-scale solutions. Randomized oversampling and forcing strategies are also tested. 
\end{abstract}

\section{Introduction}
Problems of mechanics and flow in porous media  have wide ranging applications in  many areas of science and engineering. Particularly in geomechanical modeling and its applications to reservoir engineering
for enhanced production and environmental safety due to overburden subsidence and compaction \cite{sayers2007introduction,zoback2010reservoir}. One of the key challenges is the multiscale nature of the geomechanical problems.
 Heterogeneity of reservoir properties should be accurately accounted in the geomechanical model, and this  requires a high resolution solve that adds many degrees of freedom that can be computationally costly.  Moroever, there are disparate scales between the often relatively thin reservoir structure and the large overburden  surrounding the reservoir that adds more complexity to the simulation.   
 Therefore, we propose a multiscale method  to attempt overcome some of these challenges.

The basic mathematical structure  of the poroelasticity models  are usually coupled equations for pressure and displacements known as Biot type models \cite{BiotOriginal}. 
For pressure, or flow equations, we have the  parabolic equation Darcy equation with a time dependent coupling to volumetric strain. The stress equation is the quasi-static elasticity equations with a coupling to the pressure gradients as a forcing. 
Poroelastic models of this type have been explored in the petroleum engineering  literature in the context of geomechanics for some time \cite{settari2001advances,settari1998coupled,mikelic2013convergence,minkoff2003coupled} to name just a few.
 There are noted issues that arise. The first being heterogeneities of the reservoir and surrounding media add many complications to the effective simulation due to complexity of scales. Moreover, development of flow and mechanics simulation were often considered separately. Progress was made on this problem by considering various  coupling strategies \cite{settari1998coupled}. However, in the instance that the physics is not well understood a fully coupled scheme may be desired. 
  This separation of development from flow and mechanics methods adds the complication of the computational  grids not being the same in each regime. Some effort has been made in the improvements of gridding techniques between geomechanical and flow calculations \cite{shrivastava2011use} and references therein.

As briefly noted before, typically for numerical solution of such coupled systems time splitting schemes are often used. Various splitting techniques for poroelastic equations have been explored and analyzed in the context of reservoir geomechanics in \cite{kim2011,kim2010sequential}. Also, in the context of poroelasticty and thermoelastic equations, various splitting techniques have been analyzed and implemented \cite{vab2014}. The primary splitting techniques are the undrained, fixed-stress, and fully implicit. Due to observed better errors, we will primarily consider the less computationally costly fixed-stress splitting and the more robust, yet with a loss in some matrix sparsity, fully implicit coupled approach.

Once the equations have been split in time we wish to resolve in space and will utilize a multiscale method. There are many very effective multiscale frameworks that have been developed in recent years. There are rigorous approaches based on homogenization of partial differential equations, where effective equations are derived based fine-scale equations at the microstructure level \cite{brown2011,brown2014}. However, these approaches may have limited computational use and more practical multiscale methods are used. Examples include the Heterogeneous Multiscale Method (HMM), where  macro-scale equations on coarse-grids are solved while  the effective coefficients on the fine-scale are resolved at each coarse grid nodes
 \cite{E:Engquist:2003,Abdulle:E:Engquist:Vanden-Eijnden:2012}. 
 An approach based on the Variational Multiscale Method (see \cite{MR2300286}), where coarse-grid quasi-interpolation operators are used to build an orthogonal splitting into a multiscale space and a fine-scale space \cite{MP11}.
Fine-scale space corrections are then localized to create a computationally efficient scheme. In this paper, we will use the Generalized Multiscale Finite Element Method framework, which is a generalization of the multiscale finite element method \cite{eh09}.



To efficiently solve these splitting schemes and overcome some of  the challenges of heterogenous reservoir properties and gridding issues between mechanics and flow, we will develop a Generalized Multiscale Finite Element Method (GMsFEM) \cite{egh12}. 
Our GMsFEM has the advantage of being able to capture small scale features from the heterogeneities  into coarse-grid basis functions and offline spaces, as well as having a unified computational grids for both mechanics and flow solves.  The offline multiscale basis construction may proceed in both fluid and mechanics in parallel and both constructions are comparable. 
We proceed by first generating a coarse-grid
 and in each grid block a local static problem with varying boundary conditions is solved to construct the snapshot spaces. We then perform a dimension reduction of the snapshot space by solving auxiliary eigenvalue problems. Taking the corresponding smallest eigenpairs, and multiplying by a multiscale partition of unity we are able to construct our offline basis. In this greatly reduced dimension offline basis, the online solutions may be calculated for pressure and displacements for any viable boundary condition or forcing.

The work is organized as follows. In Section 2 we provide the mathematical background of the poroelasticity problem. We will introduce the Biot type model and highlight where the heterogeneities primarily occur. In our formulation, the computational domain will be entirely inside of the fluid filled, or reservoir, region. However, coupling to regimes of pure elasticity to model the overburden are of course possible.
In  Section 3, to outline the difficulties in full direct numerical simulation we introduce the fine-scale discretizations using coupled and splitted schemes. Once we split the porooelastic system we will be able to apply our multiscale method.   In Section 4, we present our GMsFEM algorithm and outline its construction  procedure. We will use the offline multiscale basis functions to calculate accurately   pressure and displacements, at a reduced dimension and computational cost in the online phase. Finally, numerical implementations are presented in Section 5. Using the GMsFEM,  we compare the multiscale solution to fine-scale solutions and give error estimates. We will present two different examples with varying coefficients. Additionally, we will implement and discuss different strategies with  oversampling and randomized forcings to construct the multiscale spaces. 

\section{Problem formulation} 
We denote our computational domain $\Omega\subset\mathbb{R}^d$ to be a bounded Lipschitz region.  We 
 consider linear poroelasticity problem where we wish to find a pressure $p$ and displacements $u$ satisfying 
 \begin{subequations}\label{eq:main}
\begin{eqnarray}
\label{eq:poroelas}
- \div  \sigma ( u) + \alpha \grad (p)  &= 0 \text{ in } \Omega, \\
 \alpha \frac{\partial \div  u}{\partial t} + \frac{1}{M} \frac{\partial p}{\partial t} - \div \left( \frac{k}{\nu} \grad p \right) &= f \text{ in } \Omega,
\end{eqnarray}
\end{subequations}
with initial condition for pressure 
$
p( x, 0) = p_0. 
$
We write the boundary of the domain into four sections  $\partial \Omega = \Gamma_1 \cup \Gamma_2 = \Gamma_3  \cup \Gamma_4$. We suppose the following boundary conditions on each portion
\[
 \sigma  n = 0, \quad  x \in \Gamma_1, \quad 
 u = u_1, \quad  x \in \Gamma_2, 
 \]
 and 
 \[
-\frac{k}{\nu} \frac{\partial p}{\partial n}= 0, \quad  x \in \Gamma_3, \quad  
p = p_1, \quad  x \in \Gamma_4.
\]
Here the primary sources of the heterogeneities in the physical properties arise from $ \sigma$, the stress tensor and $k$, the permeability.
We denote $M$ to be the Biot modulus, $\nu$ is the fluid viscosity, and $\alpha$ is the Biot-Willis fluid-solid coupling coefficient.  
Here,  $f$ is a source term representing injection or production processes and $n$ is the unit normal to the boundary. Body forces, such as gravity, are neglected.
In the case of a linear elastic stress-strain constitutive relation we have that the stress tensor and symmetric strain gradient  may be expressed as
\[
 \sigma(u) = 2 \mu  \varepsilon( u) + \lambda \div(u) \,  \mathcal{I},
\quad
 \varepsilon(u) = \frac{1}{2} \left( \grad  u + \grad  u^T \right),
\]
where $\mu$, $\lambda$ are Lame coefficients, $\mathcal{I}$ is the identity tensor. In the case where the media has heterogeneous material properties the coefficients $\mu$ and $\lambda$ may be highly variable. 

The above poroelasticity problem \eqref{eq:poroelas}, assuming a linear elastic stress-strain relation, can be written in operator matrix form:
\begin{eqnarray}
\label{eq:mmd}
A  u + \alpha G p = 0,  \\
\frac{d}{dt} \left( S \, p + \alpha {D}  u \right) + B p = f,
\end{eqnarray}
where
\[
A  v = -\mu \nabla^2  v - (\lambda + \mu) \grad \div  v,
\quad 
B p = -\div \left( \frac{k}{\nu} \grad p \right),
\]
and $G$ and $D$ are gradient and divergence operators and $S=\frac{1}{M}{\cal I}$.

\section{Fine-Scale Discretization}
We will now present splitting methods for the above system in the context of solving the fine-scale approximation. This will highlight the areas where we would like to utilize a multiscale method when solving in the spatial variables due to the degrees of freedom required in resolving the system. 
For approximating the  numerical solution to \eqref{eq:main} on fine-scale grid we use a standard  finite element method. We begin by  giving the corresponding variational form of the continuous problem written as 
\begin{eqnarray}
\label{eq:canu}
a( u,  v) + g(p,  v)&= 0, \quad \text{for all}   \, \,  v \in \hat{ V}, \\
d \left( \frac{d  u}{dt}, q \right)  + c \left( \frac{d p}{dt}, q \right) 
+ b(p , q)&= (f , q), \quad \text{for all}  \, \,  q \in \hat Q.
\end{eqnarray}
for $u \in V$, $p \in Q$ where 
\[
V  = \{ v \in [H^1(\Omega)]^d:  v( x) = u_1,  x \in \Gamma_2  \}, 
\quad
Q  = \{  q \in H^1(\Omega):  q(x) = p_1,  x \in \Gamma_4 \},
\]
and the test spaces  with  homogeneous boundary conditions are given by
\[
\hat{V} = \{ v \in [H^1(\Omega)]^d:  v(x) = 0,  x \in \Gamma_2  \}, 
\quad
\hat{Q} = \{ q \in H^1(\Omega): q(x) = 0,  x \in \Gamma_4 \}.
\]
Here for bilinear and linear forms we have define
\[
a( u,  v)  = \int_{\Omega}  \sigma( u) \,  v dx, 
\quad
b(p, q) = \int_{\Omega} \left( \frac{k}{\nu} \grad p, \grad q \right) dx,
\quad
c(p, q) = \int_{\Omega} \frac{1}{M} \, p \, q \, dx,  
\]
and 
\[
g(p,  v)  = \int_{\Omega}\alpha (\grad p,  v) dx,
\quad
d( u, q) = \int_{\Omega} \alpha \div  u \, q \, dx , 
\quad
(f, q) = \int_{\Omega} f  \, q \, dx.
\]
Here $\left( \cdot, \cdot \right)$ under the integrand denotes the standard inner product. In Section \ref{Numerics}, we will discretize the spaces using a fine-scale standard FEM  and denote them $V_{h}, Q_{h}$ and $\hat{V}_{h}, \hat{Q}_{h}$, $h$ being  the fine-grid size. The FEM using these spaces will serve as a reference solution for our GMsFEM outlined in Section \ref{GMsFEM}.

To solve the above system we first discretize in time. This discretization leads to several possible couplings between time-steps and the two equations of prorelasticity. We proceed by giving the coupled and so-called fixed-stress splitting \cite{kim2011,vab2014}.
The standard fully implicit finite-difference scheme, or coupled scheme, can be used for the time-discretization and is given by 
\begin{subequations}\label{coupled}
\begin{eqnarray}
\label{eq:capp}
a( u^{n+1},  v) + g(p^{n+1},  v) = 0,  \\
d \left( \frac{ u^{n+1} -  u^n}{\tau}, q \right)  
+ c \left( \frac{p^{n+1} - p^n}{\tau}, q \right)
+ b(p^{n+1} , q) = (f , q),  
\end{eqnarray}
\end{subequations}
with $u^n =  u(x, t_n)$,  $p^n = p(x, t_n)$, where $t_n = n \tau$, $n = 0, 1, ..., M_{T}$, $M_{T} \tau = T$ and $\tau > 0$.
%
For time discretization we can apply many different splitting techniques which often occur in the literature.  

Another we shall consider here is the fixed-stress splitting scheme
\begin{subequations}\label{fixedstress}
\begin{eqnarray}
\label{eq:sapp}
a( u^{n+1},  v) + g(p^{n+1},  v) = 0,  \\
\label{pressuresplit}
d \left( \frac{ u^{n} - u^{n-1}}{\tau}, q \right)  
+ s \left( \frac{p^{n+1} - p^n}{\tau}, q \right)
+ b(p^{n+1} , q) = l(q), 
\end{eqnarray}
\end{subequations}
where the variational form is re-written with
\[
s(p, q) = \int_{\Omega} \left( \frac{1}{M} + \frac{\alpha^2}{K_{dr}} \right) \, p \, q \, dx,  
\quad
l(q) = \int_{\Omega} 
\left( f  + \frac{\alpha^2}{K_{dr}} \frac{p^n - p^{n-1}}{\tau}  \right) \, q \, dx
\]
and $K_{dr}$ is the drained modulus
\[
K_{dr} = \frac{E (1- \nu_p)}{(1-2 \nu_p) (1+\nu_p)},
\]
where $\nu_p$ is the Poisson ratio and $E$ is the elastic modulus. When we utlize the fixed-stress splitting scheme, first we solve pressure equation for $p^{n+1}$ given data at the previous time-steps. Then, passing this new pressure information, we return to the quasi-static stress equation and  calculate displacements at $u^{n+1}$.

\section{GMsFEM for Poroelasticity} \label{GMsFEM}
In the GMsFEM presented here, we will focus on the development in the fixed-stress splitting \eqref{fixedstress}. We will however give numerical examples from both coupling strategies.
The fixed-stress splitting decouples the flow and mechanics equations. We will first present the offline multiscale basis construction in the fluid or pressure solve then its construction in the mechanics or displacement calculation step. In this algorithm, due to the heterogeneities arising primarily from the permeability $k$ and the stress tensor $\sigma(u)$, we will solve local problems in each of the relevant portions of the variational form to construct the offline multiscale spaces. 
We now outline the general procedure of the GMsFEM algorithm.

The overall fine-scale model equations will be solved on a fine-grid using spaces $V_{h}, Q_{h}$ and $\hat{V}_{h}, \hat{Q}_{h}$, and will act as our reference solutions. Once the fine-grid is established we must introduce the concepts of coarse-grids and their relationships.
%
%
To this end, let $\mathcal{T}^H$ be a standard conforming partition of the computational domain $\Omega$ into finite elements. 
We refer to this partition as the coarse-grid and assume that each coarse element is partitioned into a connected union of fine grid blocks. 
The fine grid partition will be denoted by $\mathcal{T}^h$, and is by definition a refinement of the coarse grid $\mathcal{T}^H$.
We use $\{x_i\}_{i=1}^{N}$, where $N$ is the number of coarse nodes, to denote the vertices of the coarse mesh $\mathcal{T}^H$, and define the neighborhood of the node $x_i$ by
\[
\omega_i=\bigcup_{j}\left\{ K_j\in\mathcal{T}^H \, | \,  x_i\in \overline{K}_j\right\}.
\]
\begin{figure}[htb]
  \centering
  \includegraphics[width=0.6\textwidth]{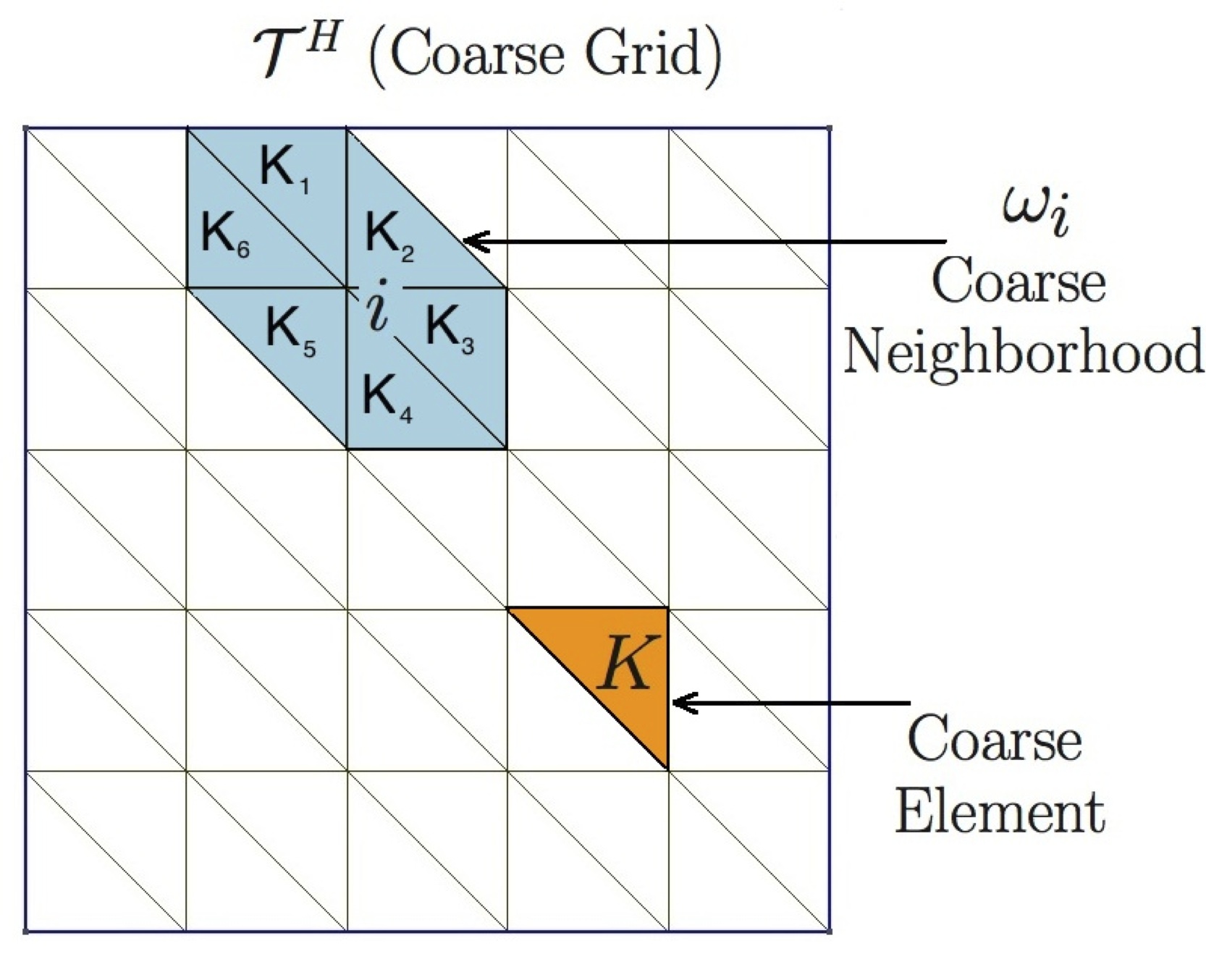}
  \caption{Illustration of a coarse neighborhood and coarse element}
  \label{schematic}
\end{figure}
See Figure~\ref{schematic} for an illustration of neighborhoods and elements subordinated to the coarse discretization. We emphasize that the use of $\omega_i$ is to denote a coarse neighborhood, and we use $K$ to denote a coarse element throughout the paper.

Boadly speaking, the GMsFEM algorithm consist of several steps:
\begin{itemize}
\item {\bf Step 1: } Generate the coarse-grid, $\mathcal{T}^H$.
\item {\bf Step 2: } Construct the  snapshot space, used to compute an offline space, by solving many local problems on the fine-grid.
\item {\bf Step 3: } Construct  a small dimensional offline space by performing dimension reduction in the space of local snapshots.
\item {\bf Step 4: } Use small dimensional offline space to find the solution of a coarse-grid problem for any force term and/or boundary condition.
\end{itemize}

As noted previously, because coupled system of equations for poroelasticity problems can be solved using splitting scheme, we can construct multisclate basis functions for pressure and displacements separately. We begin by considering the pressure solve, then, the displacement solve. 

\subsection{Pressure Solve}\label{pressuresolve}
Recall, for the  numerical solution of pressure equation on coarse grid we consider the continuous Galerkin (CG) formulation \eqref{pressuresplit} given by
\begin{equation} 
\label{eq:p-discr}
s( \frac{p^{n+1} - p^n}{\tau}, q) +b(p^{n+1}, q) = l(q) - d \left( \frac{ u^{n} - u^{n-1}}{\tau}, q \right)  , \quad 
\text{for all} \, \, q \in Q_{\text{off}},
\end{equation}
where  $Q_{\text{off}}$ is used to denote the space spanned by multiscale basis functions $\psi_k^{\omega_i}$, each of which is supported in $\omega_i$.
The index $k$ represents the numbering of these multiscale  basis functions. We will now show how to construct the offline multiscale space $Q_{\text{off}}$.
In turn, the CG solution of the form  
\[
p(x, t)=\sum_{i,k} p_{k}^i(t) \psi_{k}^{\omega_i}(x),
\]
will be sought.

We begin by construction of a snapshot space $V_{\text{snap}}^{\omega}$.  We use harmonic extensions
\begin{equation} 
\label{harmonic_ex}
\begin{split}
b(\psi_{j}^{\omega, \text{snap}} , q)&= 0 \quad \text{ in } \, \omega, \\
\psi_{j}^{\omega, \text{snap}}&=\delta_j^h(x) \quad \text{ on } \partial\omega.
\end{split}
\end{equation}
Here $\delta_j^h(x)$  are defined by
$\delta_j^h(x)=\delta_{j,k},\,\forall j,k\in \textsl{J}_{h}(\omega_i)$, where $\textsl{J}_{h}(\omega_i)$ denotes the fine-grid boundary node on $\partial\omega_i$. For simplicity, we will omit the index $i$ when there is no ambiguity.

Let $l_i$ be the number of functions in the snapshot space in the region $\omega$, and define
\[
Q^{\omega}_{\text{snap}} = \text{span}\{ \psi_{j}^{ \text{snap}}:\quad 1\leq j \leq l_i \},
\]
for each coarse subdomain $\omega$.
We denote the corresponding matrix of snapshot functions to be 
$$
R^p_{\text{snap}} = \left[ \psi_{1}^{\text{snap}}, \ldots, \psi_{l_i}^{\text{snap}} \right].
$$

To construct the offline space $Q_{\text{off}}$, we perform a dimension reduction of the space of snapshots by using an auxiliary spectral decomposition. More precisely,  we solve the  eigenvalue problem in the space of snapshots:
\begin{equation} 
\label{offeig}
B^{\text{off}} \Psi_k^{\text{off}} = \lambda_k^{\text{off}} M^{\text{off}} \Psi_k^{\text{off}},
\end{equation}
where
\[
B^{\text{off}} = (R^p_{\text{snap}})^T B R^p_{\text{snap}},
\quad
M^{\text{off}} =  (R^p_{\text{snap}})^T M R^p_{\text{snap}},
\]
where $B$ and $M$ denote fine scale matrices
\[
B_{ij} = \int_{\Omega} \left( \frac{k}{\nu} \grad \phi_i , \grad \phi_j \right)\, dx, 
\quad
M_{ij} = \int_{\Omega}  \frac{k}{\nu}  \phi_i \phi_j \, dx.
\]
Here, $\phi_i$ are  fine-scale basis functions.

We then  choose the smallest $N^{\omega,p}_{\text{off}}$ eigenvalues from Eq.~\eqref{offeig} and form the corresponding eigenvectors in the space of snapshots by setting
$$\psi_k^{\text{off}} = \sum_{j=1}^{l_i} \Psi_{kj}^{\text{off}} \psi_j^{\text{snap}},$$
 for $k=1,\ldots, N^{\omega,p}_{\text{off}}$, where $\Psi_{kj}^{\text{off}}$ are the coordinates of the vector $\psi_{k}^{\text{off}}$. We denote the span of this reduced space as $Q_{\text{off}}^{\omega}$.
 
%
For construction of the offline space, to ensure the functions we construct form an $H^1$ conforming basis, we define multiscale partition of unity functions $\chi_i$ 
\begin{eqnarray} 
\label{pou}
b(\chi_i , q) = 0 \quad \text{ in } K, \\
\chi_i = g_i \quad \text{on} \, \, \, \partial K, \nonumber
\end{eqnarray}
for all $K \in \omega$. Here $g_i$ is a continuous on $K$ and is linear on each edge of $\partial K$. We could choose $g_i$ to also be selected shape function, Neumann conditions, or boundary conditions on larger domains in the context of oversampling. 
%
 
 
Finally, we  multiply the partition of unity functions by the eigenfunctions in the offline space $Q_{\text{off}}^{\omega_i}$ to construct the resulting basis functions
\begin{equation} 
\label{cgbasis}
\psi_{i,k} = \chi_i \psi_k^{\omega_i, \text{off}} \quad \text{for} \, \, \,
1 \leq i \leq N \, \, \,  \text{and} \, \, \, 1 \leq k \leq N_{\text{off}}^{\omega_i,p},
\end{equation}
where $N_{\text{off}}^{\omega_i,p}$ denotes the number of offline eigenvectors that are chosen for each coarse node $i$. We note that the construction in Eq.~\eqref{cgbasis} yields  continuous basis functions due to the multiplication of offline eigenvectors with the initial (continuous) partition of unity. Next, we define the continuous Galerkin spectral multiscale space as
\begin{equation} \label{cgspace}
Q_{\text{off}}  = \text{span} \{ \psi_{i,k} : \,  \, 1 \leq i \leq N \, \, \,  \text{and} \, \, \, 1 \leq k \leq N_{\text{off}}^{\omega_i,p}  \}.
\end{equation}
Using a single index notation, we may write $Q_{\text{off}} = \text{span} \{ \psi_{i} \}_{i=1}^{N^p_c}$, where $N^p_c =\sum_{i=1}^{N}N_{\text{off}}^{\omega_{i},p}$
denotes the total number of basis functions in the spaces $Q^{\omega_{i}}_{\text{off}}$, for  $i=1,\dots,N$. 

Denote the matrix  
\[
R^p = \left[ \psi_1 , \ldots, \psi_{N^p_c} \right]^T,
\] 
where $\psi_i$ are used to denote the nodal values of each basis function defined on the fine grid.
Then, the variational form in \eqref{eq:p-discr} yields the following linear algebraic system
\begin{equation}\label{coarse.pressure}
Q_c p_c^{n+1} = Y_c p_c^{n},
\end{equation}
where 
$$Q_c = R^p (\frac{1}{M \tau} + B ) (R^p)^T , \quad Y_c = R^p F_p.$$ Here,  $F_{p}$ being the operator corresponding right hand side data from the previous time step and $p_c$ denotes the coarse-scale nodal values of the discrete CG solution.     
We also note that the operator matrix may be analogously used in order to project coarse scale solutions onto the fine grid $$p^{n+1} =  (R^p)^T p_c^{n+1}.$$

\subsection{Displacement Solve}
We now suppose that we have solved for the fine-grid pressure $p^{n+1}$ by the GMsFEM pressure solve in the previous section. We must now solve the mechanics equations \eqref{eq:sapp}. Since the construction of the multiscale offline space remains very similar in this setting, we will be a bit more brief on its construction. 
Recall, for discretization of the displacements equation we rewrite equation as follows
\begin{equation}
\label{eq:ell}
{A}  u^{n+1}   = F_u, 
\end{equation}
where $F_u = -\alpha {G} p^{n+1}$. The corresponding continuous Galerkin (CG) formulation for displacements equations is given by:
\begin{equation} 
\label{eq:u-discr}
a(u^{n+1}, v) = (f_u, v), \quad 
\text{for all} \, \, v\in V_{\text{off}},
\end{equation}
where  $u(x, t)=\sum_{i,k} u_{k}^i(t) \varphi_{k}^{\omega_i}(x)$,  where $\varphi_{k}^{\omega_i}$ are fine-scale basis functions, and we  construct the multiscale offline space $V_{\text{off}}$.

For construction of multiscale basis functions for displacements we use similar algorithm that we used for pressure. 
For construction of a snapshot space $V_{\text{snap}}^{\omega}$  we solve following problem in $\omega$
\begin{equation} 
\label{harmonic_ex2}
\begin{split}
a(\varphi_{j}^{\omega, \text{snap}} , v) &= 0 \quad \text{in } \, \omega, \\
\varphi_{j}^{\omega, \text{snap}}&=\delta_j^h(x), \quad \text{ on }\partial\omega.
\end{split}
\end{equation}
Let $l_i$ be the number of functions in the snapshot space in the region $\omega$, and define
\[
V^{\omega}_{\text{snap}} = \text{span}\{ \varphi_{j}^{ \text{snap}}:\quad 1\leq j \leq l_i \},
\]
for each coarse subdomain $\omega$. Note we are using the same notation but with different harmonic extensions.
We denote the corresponding matrix of snapshot functions, again with similar notation, to be 
$$
R^u_{\text{snap}} = \left[ \varphi_{1}^{\text{snap}}, \ldots, \varphi_{l_i}^{\text{snap}} \right].
$$

Again, we perform a dimension reduction of the space of snapshots by using an auxiliary spectral decomposition. We solve the  eigenvalue problem in the space of snapshots
\begin{equation} 
\label{offeig2}
A^{\text{off}} \Phi_k^{\text{off}} = \lambda_k^{\text{off}} N^{\text{off}} \Phi_k^{\text{off}},
\end{equation}
where
where
\[
A^{\text{off}} = (R^u_{\text{snap}})^T A R^u_{\text{snap}},
\quad
N^{\text{off}} = ( R^u_{\text{snap}})^T N R^u_{\text{snap}},
\]
where $A$ and $N$ denote fine scale matrices
\[
A_{mn}
= \int_{\Omega} \Big( 2\mu  {\varepsilon}( \varphi_m^{\text{}}) :  {\varepsilon}( \varphi_n^{\text{}})
+ \lambda \div  ({\varphi_m^{\text{}}}) \cdot\div(  \varphi_n^{\text{}} ) \Big),
\]
and
\[
 N_{mn}
= \int_{\Omega}  (\lambda + 2 \mu)  \varphi_m^{\text{}} \cdot  \varphi_n^{\text{}}.\]
Here, $\varphi_i$ are  fine-scale basis functions.

We then  choose the smallest $N^{\omega,u}_{\text{off}}$ eigenvalues from Eq.~\eqref{offeig2} and form the corresponding eigenvectors in the space of snapshots by setting
$$\varphi_k^{\text{off}} = \sum_{j=1}^{l_i} \Phi_{kj}^{\text{off}} \varphi_j^{\text{snap}},$$
 for $k=1,\ldots, N^{\omega,u}_{\text{off}}$, where $\Phi_{kj}^{\text{off}}$ are the coordinates of the vector $\varphi_{k}^{\text{off}}$. We denote the span of this reduced space as $V_{\text{off}}^{\omega}$.

For construciton of multiscale partition of unity functions for the mechanics solve, we proceed as before and solve for all $K \in \omega$
\begin{eqnarray} 
\label{pou2}
a(\xi_i , v) = 0 \quad \text{ in } K,  \\
\xi_i = g_i \quad \text{on} \, \, \, \partial K.
\end{eqnarray}
Here $g_i$ is a continuous function on $K$ and is linear on each edge of $\partial K$.
Finally, we  multiply the partition of unity functions by the eigenfunctions in the offline space $V_{\text{off}}^{\omega_i}$ to construct the resulting basis functions
\begin{equation} 
\label{cgbasis.mechanics}
\varphi_{i,k} = \xi_i \varphi_k^{\omega_i, \text{off}} \quad \text{for} \, \, \,
1 \leq i \leq N \, \, \,  \text{and} \, \, \, 1 \leq k \leq N_{\text{off}}^{\omega_i,u},
\end{equation}
where $N_{\text{off}}^{\omega_i,u}$ denotes the number of offline eigenvectors that are chosen for each coarse node $i$. 
Next, we define the spectral multiscale space as
\begin{equation} \label{cgspace.mechanics}
V_{\text{off}}  = \text{span} \{ \varphi_{i,k} : \,  \, 1 \leq i \leq N \, \, \,  \text{and} \, \, \, 1 \leq k \leq N_{\text{off}}^{\omega_i}  \}.
\end{equation}
Using a single index notation, we may write $V_{\text{off}} = \text{span} \{ \varphi_{i} \}_{i=1}^{N^u_c}$, where $N^u_c =\sum_{i=1}^{N}N_{\text{off}}^{\omega_{i},u}$
denotes the total number of basis functions in the space $V^{\omega_{i}}_{\text{off}}$, for all $i=1,\dots,N$.

And after construction $V_{\text{off}}$ we 
denote the matrix  
\[
R^u = \left[ \varphi_1 , \ldots, \varphi_{N^u_c} \right]^T,
\] 
where $\varphi_i$ are used to denote the nodal values of each basis function defined on the fine grid.
Then, the variational form in \eqref{eq:u-discr} yields the following linear algebraic system
\begin{equation}\label{coarse.mechanics}
A_c u_c^{n+1} = F_c,
\end{equation}
where $A_c = R^u A (R^u)^T$, $F_c = R^u F_u$ and $$u^{n+1} =  (R^u)^T u_c^{n+1}.$$

\section{Numerical Examples}\label{Numerics}
In this section, we present numerical examples to  demonstrate the  performance of the GMsFEM  for computing the solution of the poroelasticity problem in heterogenous domains. Although we presented the algorithm in the fixed-stress splitting, we are able to apply the same offline spaces $(Q_{\text{off}},V_{\text{off}})$ as their construction remains the same in the fully coupled setting.  However, in the coupled setting the equations \eqref{coarse.pressure} and \eqref{coarse.mechanics} will no longer be decoupled and must be solved simultaneously.

We will implement a single complicated geometry with contrasting parameter values. We provide two cases one with lower contrast in elastic properties and another with higher contrast.
We present the algorithm applied to these heterogenous coefficients in both the fully coupled and fixed stress time splittings. We give the errors with varying multiscale basis functions and over time. 
We then will apply the GMsFEM method with oversampling and with  snapshots with randomized boundary conditions to obtain good accuracy, while having to solve fewer snapshot solutions. The effects of higher contrast in properties will also be discussed. 
\subsection{GMsFEM Implementation}
\begin{figure}[htb]
\begin{center}
\includegraphics[width=0.35\linewidth]{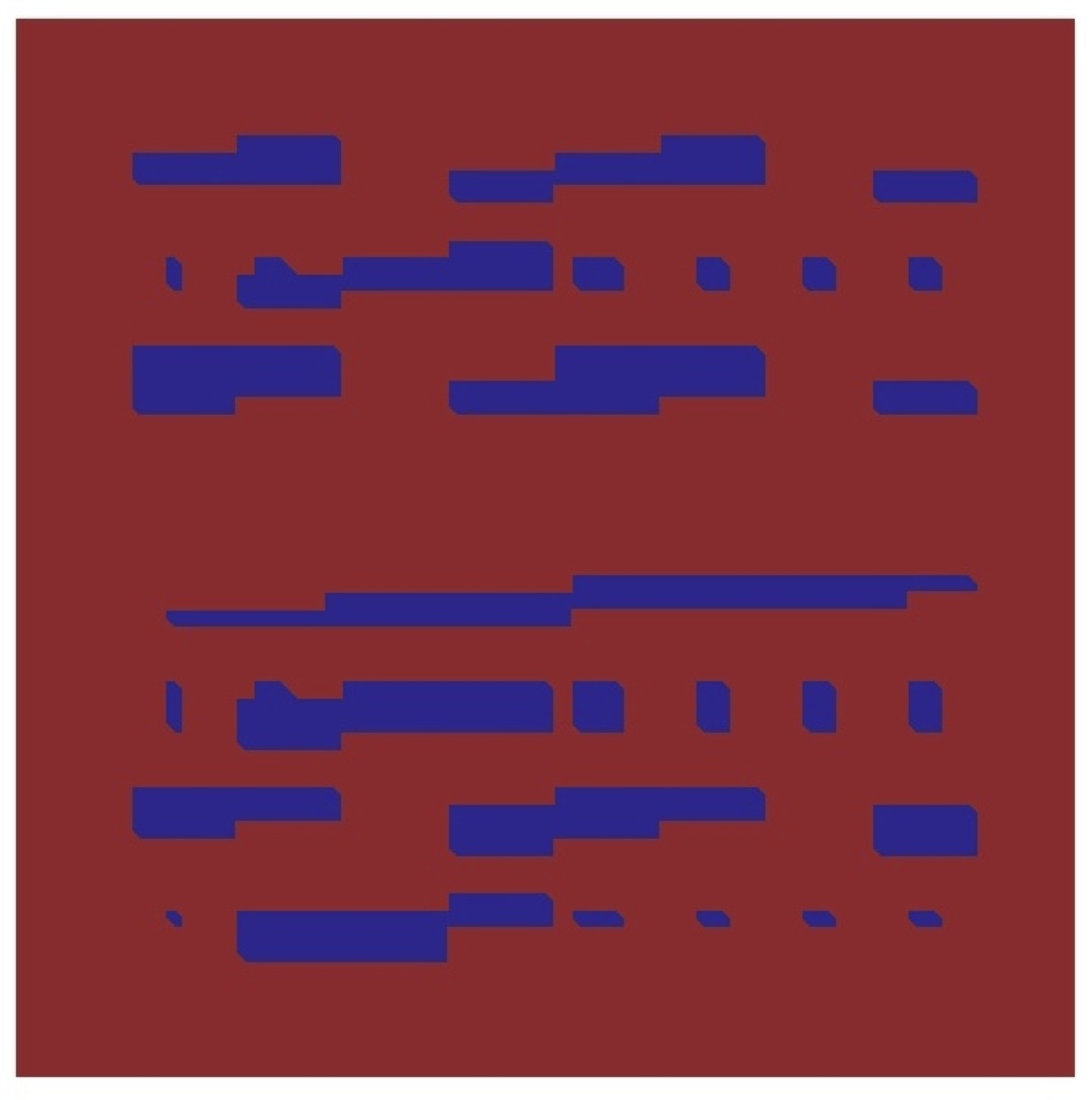} 
\caption{Coefficients subdomains. Red is the subdomain 1 and blue is the subdomain 2}
\label{fig:koeff}
\end{center}
\end{figure}

First, we take the computational domain $\Omega$ as a unit square $[0,1]^2$, and set the source term $f = 0$ in \eqref{eq:main}.
We utilize  heterogeneous coefficients that have different values in  two subdomains.  We denote each region as subdomain 1 and 2,
and use following coefficients: for the Biot modulus we take $M_1 = 1.0, M_2 = 10$ and for permeability  $k_1=10^{-3}, k_2 = 1$ in the two regions.
For fluid viscosity we take $\nu = 1$ and  fluid-solid coupling constant $\alpha = 0.9$. 
For the elastic properties, we present results for two test cases. In Case 1, the elastic modulus is given by  $E_1 = 10, E_2 = 1$ in each respective subdomain and in Case 2,  we have $E_1 = 10, E_2 = 10^{-3}$.
The Poisson's ratio is $\eta = 0.22$, and these can be related to the parameters $\mu_{i}$ and $\lambda_{i}$, for $i=1,2,$ via the relation
\[
\mu_{i} = \frac{E_{i}}{2 (1 + \eta)}, \quad
\lambda_{i} = \frac{E_{i} \eta}{(1+ \eta) ( 1- 2 \eta)},
\]
in each subdomain. The subdomains for coefficients shown in Fig. \ref{fig:koeff},  where the background media in red is the subdomain 1, and  isolated particles  and strips in blue are the subdomain 2.

As we have chosen $f=0$ we must use boundary conditions to force flow and mechanics. In these tests, we use following boundary conditions: 
\[
p = p_1, \quad x \in \Gamma_T, \quad
p = p_0, \quad x \in \Gamma_B, \quad
\frac{\partial p}{\partial n} = 0, \quad x \in \Gamma_L \cup \Gamma_R,
\] 
and
\[
u_{x} = 0, \quad \frac{\partial u_{y}}{\partial y} = 0, \quad x \in \Gamma_L, \quad
\frac{\partial u_x}{\partial x} = 0, \quad u_{y} = 0, \quad x \in \Gamma_B, 
\]
and finally,
\[
\frac{\partial u_x}{\partial x} = 0, \quad \frac{\partial u_y}{\partial y} = 0, \quad x \in \Gamma_T \cup \Gamma_R.
\]
Here $\Gamma_L$ and $\Gamma_R$ are left and right boundaries, $\Gamma_T$ and $\Gamma_B$ are top and bottom boundaries respectively. We set $p_0 = 0$ and $p_1 = 1$ to drive the flow, and thus, the mechanics. 

\begin{figure}[htb]
\begin{center}
\includegraphics[width=0.8\linewidth]{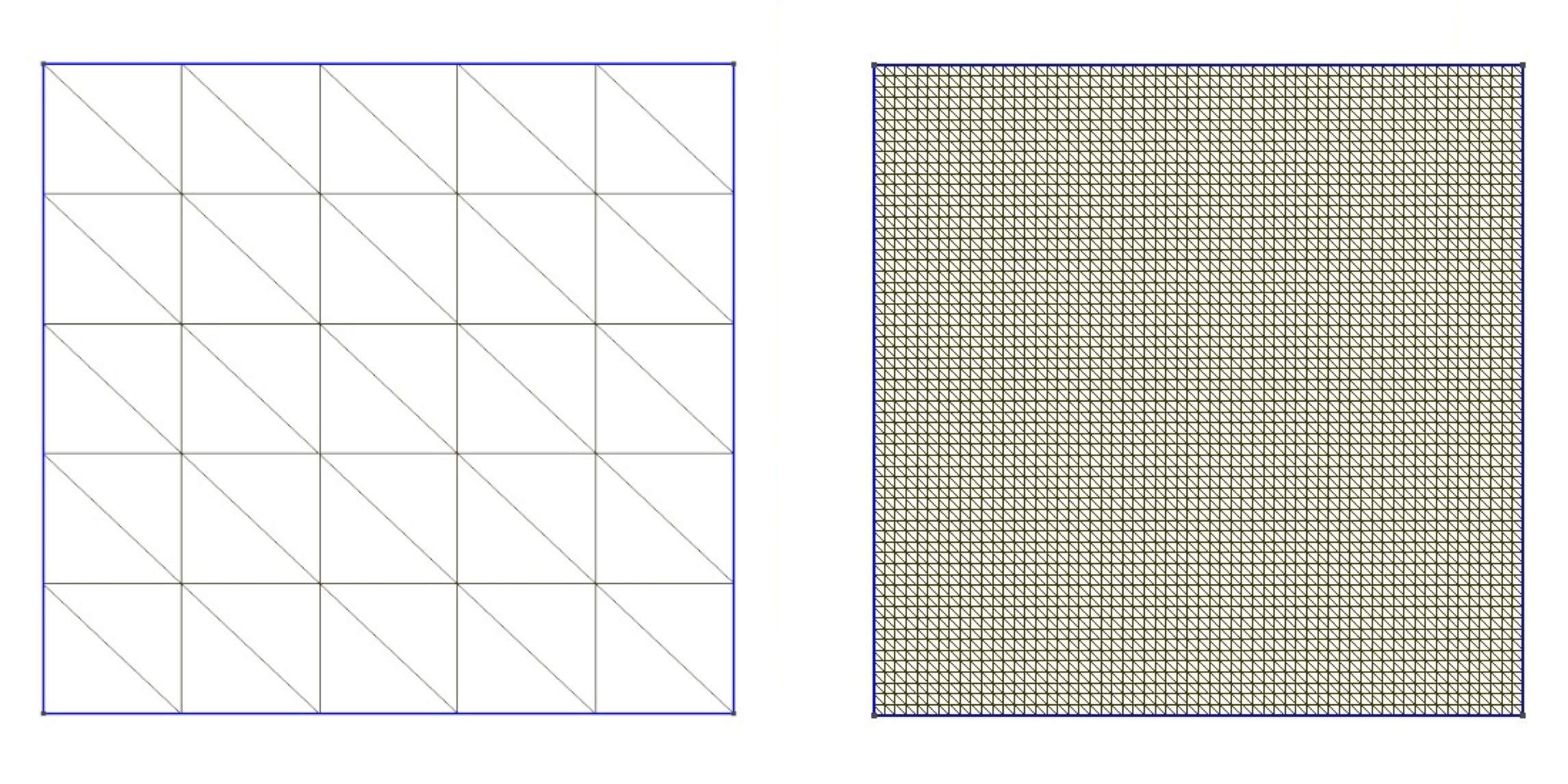} 
\caption{Coarse and fine grids}
\label{fig:domain}
\end{center}
\end{figure}

\begin{figure}[htb]
\begin{center}
\includegraphics[width=0.8\linewidth]{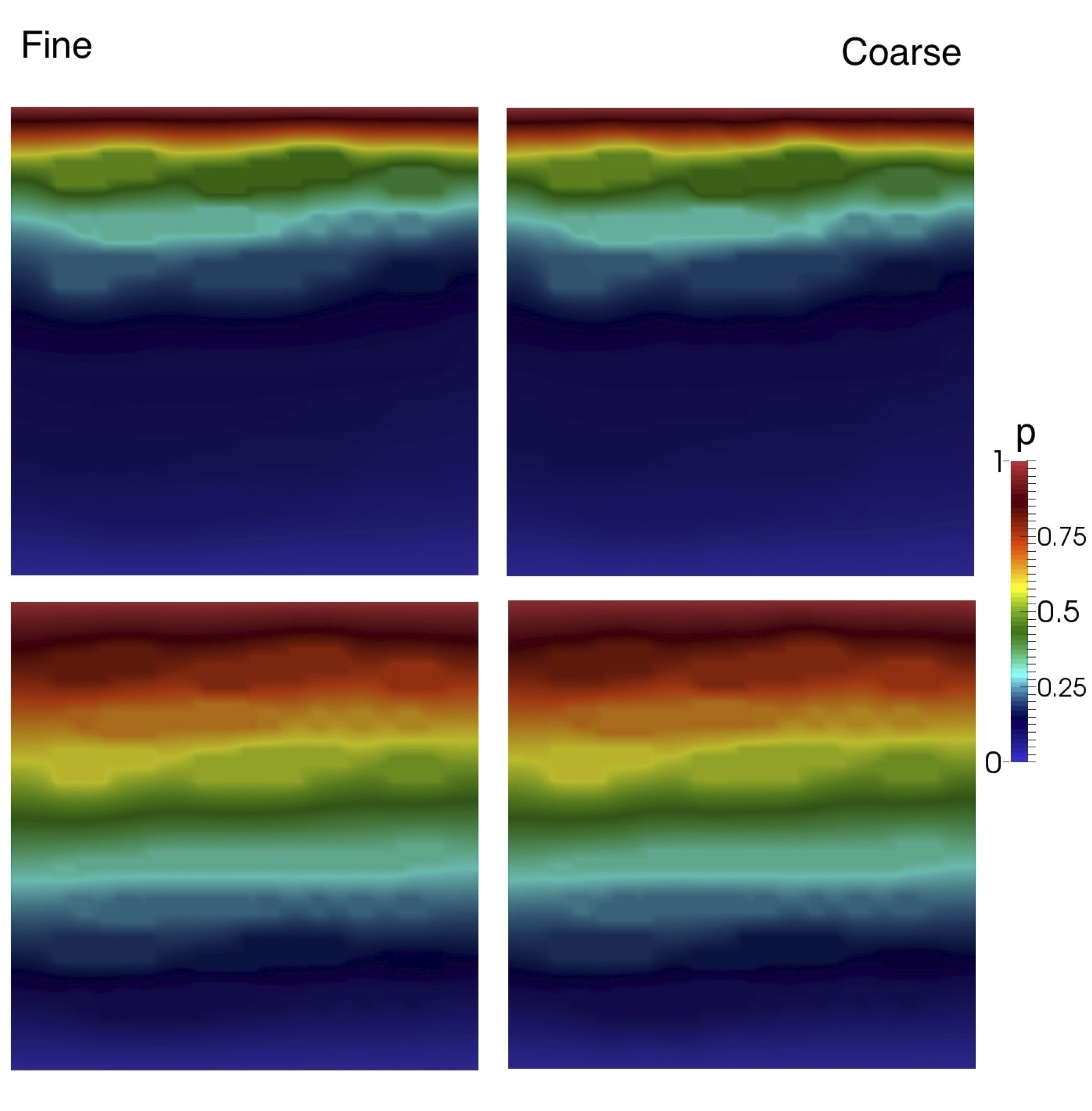} 
\caption{The fine-scale and coarse-scale solutions of the pressure distribution for $T = 10$ and $100$ (from top to bottom) for case 1. The dimension of the fine-scale solution is 11163 and the dimension of the coarse space is 864.}
\label{fig:p}
\end{center}
\end{figure}

\begin{figure}[htb]
\begin{center}
\includegraphics[width=0.8\linewidth]{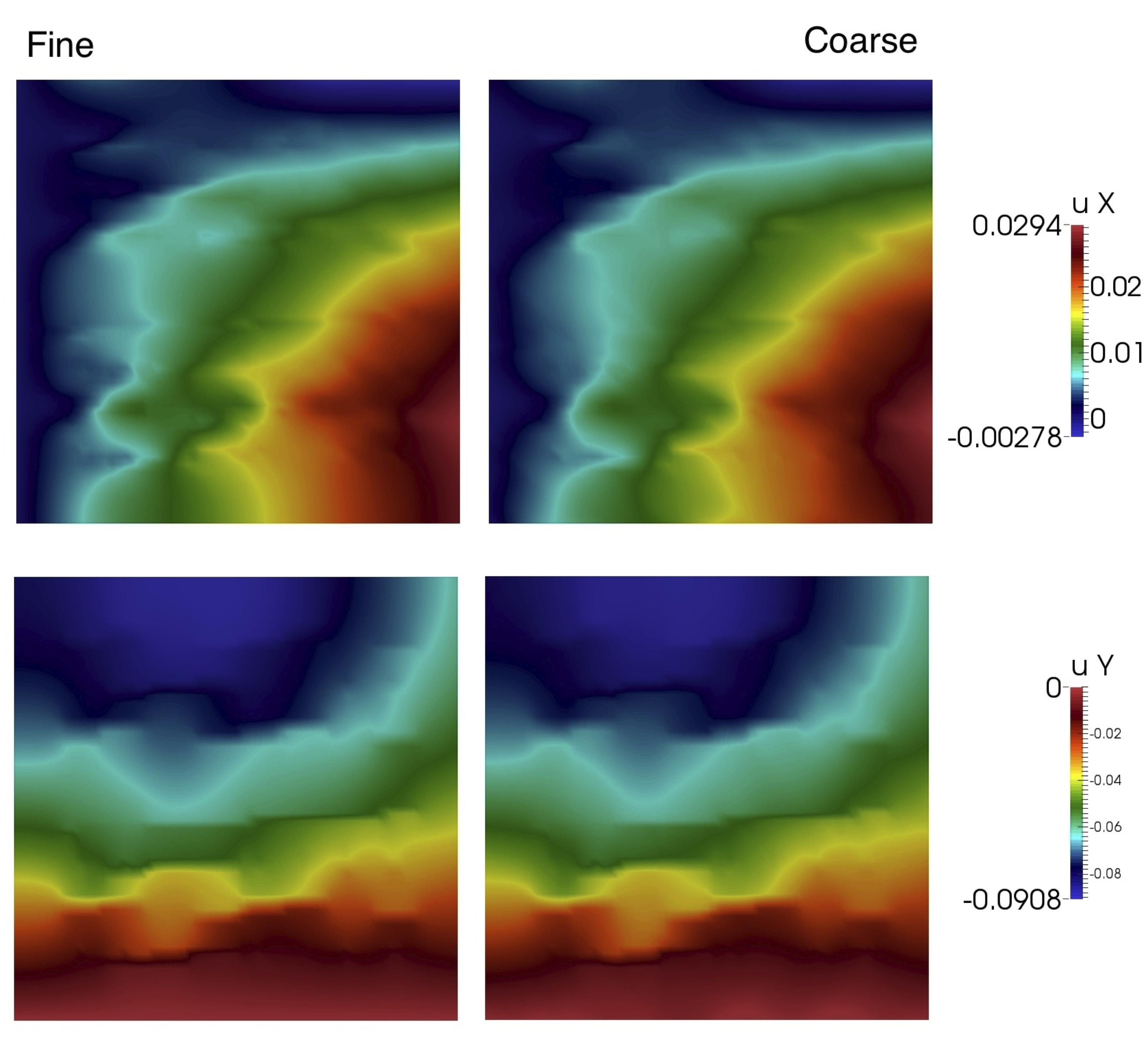} 
\caption{The fine-scale and coarse-scale solutions of the displacements $u_x$ and $u_y$ for case 1. The dimension of the fine-scale solution is 11163 and the dimension of the coarse space is 864.}
\label{fig:u}
\end{center}
\end{figure}

In Fig. \ref{fig:domain} we show the coarse and fine grids. The coarse grid consists of 36 nodes and 50 triangle cells, and the fine mesh consists of 3721 nodes and 7200 triangle cells. 
The number of time steps  is $20$ and the maximal time being set at $T_{max} = 100$. As an  initial condition for pressure we use $p = p_0$.  
The reference solution computed by using a  standard FEM (linear basis functions for pressure and displacements) on the fine grid and  using a fully coupled scheme.
The pressure and the displacement fields for  Case $1$ on the fine-grid are presented on the left column of  Fig. \ref{fig:p} - \ref{fig:u}.

We test the fully coupled and  fixed-stress splitting schemes. The errors will be measured in weighted $L^2$ and weighted $H^1$ norm and semi-norm  for pressure
\begin{align*}
\norm{e_p}_{L^2(\Omega)} &= \left(  \int_{\Omega} \frac{k}{\nu} (p_f - p_{ms})^2 dx \right)^{1/2}, \\
\seminorm{e_p}_{H^1(\Omega)} &= \left(  \int_{\Omega} \left( \frac{k}{\nu} \grad (p_f - p_{ms}), \grad (p_f - p_{ms}) \right) dx \right)^{1/2},
\end{align*}
and for displacements
\begin{align*}
\norm{e_u}_{L^2(\Omega)} &= \left(  \int_{\Omega} (\lambda + 2 \mu) (u_f - u_{ms}, u_f - u_{ms}) dx \right)^{1/2}, \\
\seminorm{e_u}_{H^1(\Omega)} &= \left(  \int_{\Omega} \left( \sigma (u_f - u_{ms}), \varepsilon(u_f - u_{ms}) \right) dx  \right)^{1/2}.
\end{align*}
Here $(u_f,p_f)$ and $(u_{ms},p_{ms})$ are  fine-scale and coarse-scale using GMsFEM solutions, respectively for pressure and displacements.

Recall, we will use a few multiscale basis functions per each coarse node $
\omega_{i}$, and these number of coarse basis defines the problem size (dimension of offline spaces, $Q_{\text{off}}$ and $V_{\text{off}}$).  
We suppose that in each patch $\omega_{i}$ we take the same number of multiscale basis functions for pressure, $N^p_{\text{off}}=N^{\omega_{i},p}_{\text{off}}$, for $i=1,\cdots,N$.
Similarly for displacements we take  $N^u_{\text{off}}=N^{\omega_{i},u}_{\text{off}}$, for $i=1,\cdots,N$.
Varying the basis functions in both pressure and displacement multiscale spaces we recorded the errors at the final times. 

In Tables \ref{tab:c1-c} and  \ref{tab:c2-c}, we present the  weighted $L^2$ and $H^1$ errors for Case 1 and Case 2 of the coefficients in geometry Fig. \ref{fig:koeff}
using the fully coupled scheme. We compare these to a fine-scale solution space with dimension 11163.
%
In these tables, $N_{\text{off}}^{p}$ and $N_{\text{off}}^{u}$ are number of multiscale basis functions for each neighborhoods, the second column show the dimension of the offline space, the next two columns present the weighted $L^2$ and $H^1$ errors for pressure and last two columns show the weighted $L^2$ and $H^1$ errors for displacements.
We see that the errors in pressure remain similar in both cases because the permeability parameters remain the same and  the change is in elastic properties between scenarios.
 In Case 2, pictured in Table \ref{tab:c2-c}, we see great errors in displacements throughout when compared to Case 1 in Table \ref{tab:c1-c} because the elastic properties in Case 2 have several orders of higher contrast.
 
In a similar setting, we consider the fixed-stress splitting. For Case 1 we present the results in  Table \ref{tab:c1-s}, the errors are very similar compared to the corresponding fully coupled scheme. This may be because we are comparing a fine-scale fully coupled scheme to a multiscale fully coupled scheme  and similarly, a fine-scale splitting scheme to a multiscale splitting scheme and the errors do not differ very much between the two schemes here. 
  For Case 2 we present the errors in Table \ref{tab:c2-s} and again see that the errors are higher when compared to the lower contrast scenario. Comparing these results with the Case 2 using the fully coupled scheme, presented in Table  \ref{tab:c2-c}, we see that both the pressure errors and displacement errors are much greater in this sequential coupling. This disparity is  particularly striking when few multiscale basis functions are used. 

\begin{figure}
\begin{center}
\includegraphics[width=0.45\linewidth]{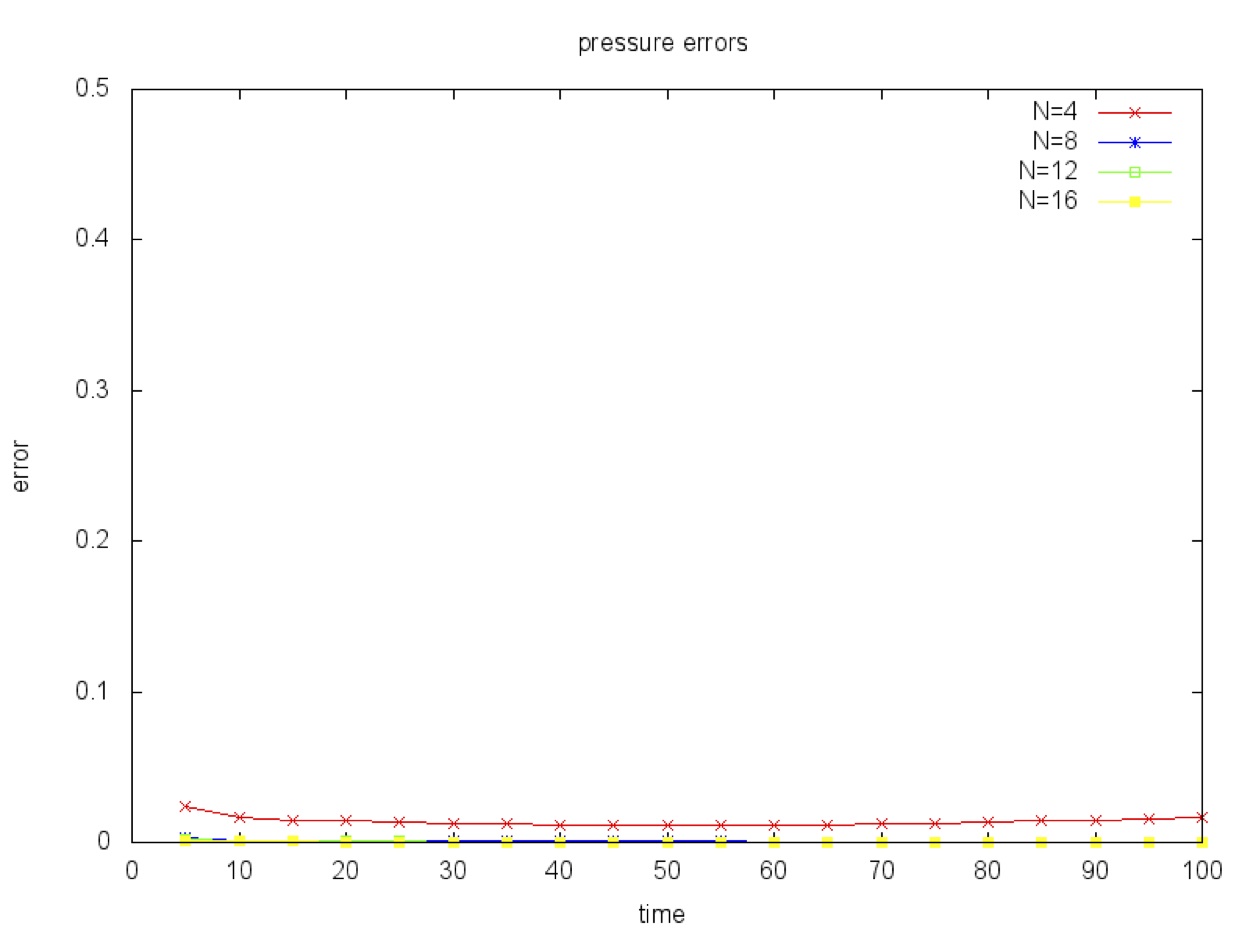} 
\includegraphics[width=0.45\linewidth]{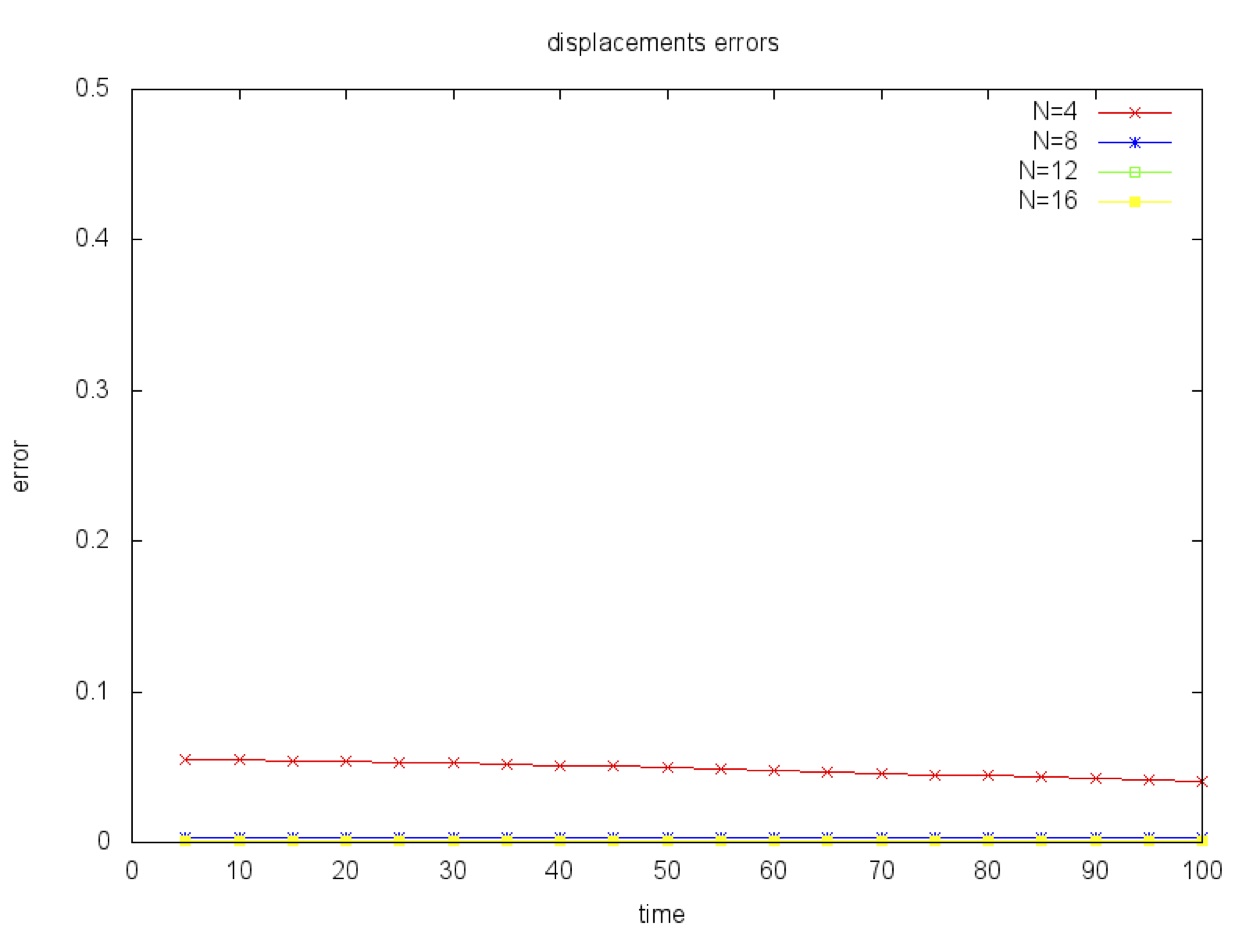} 
\\
\includegraphics[width=0.45\linewidth]{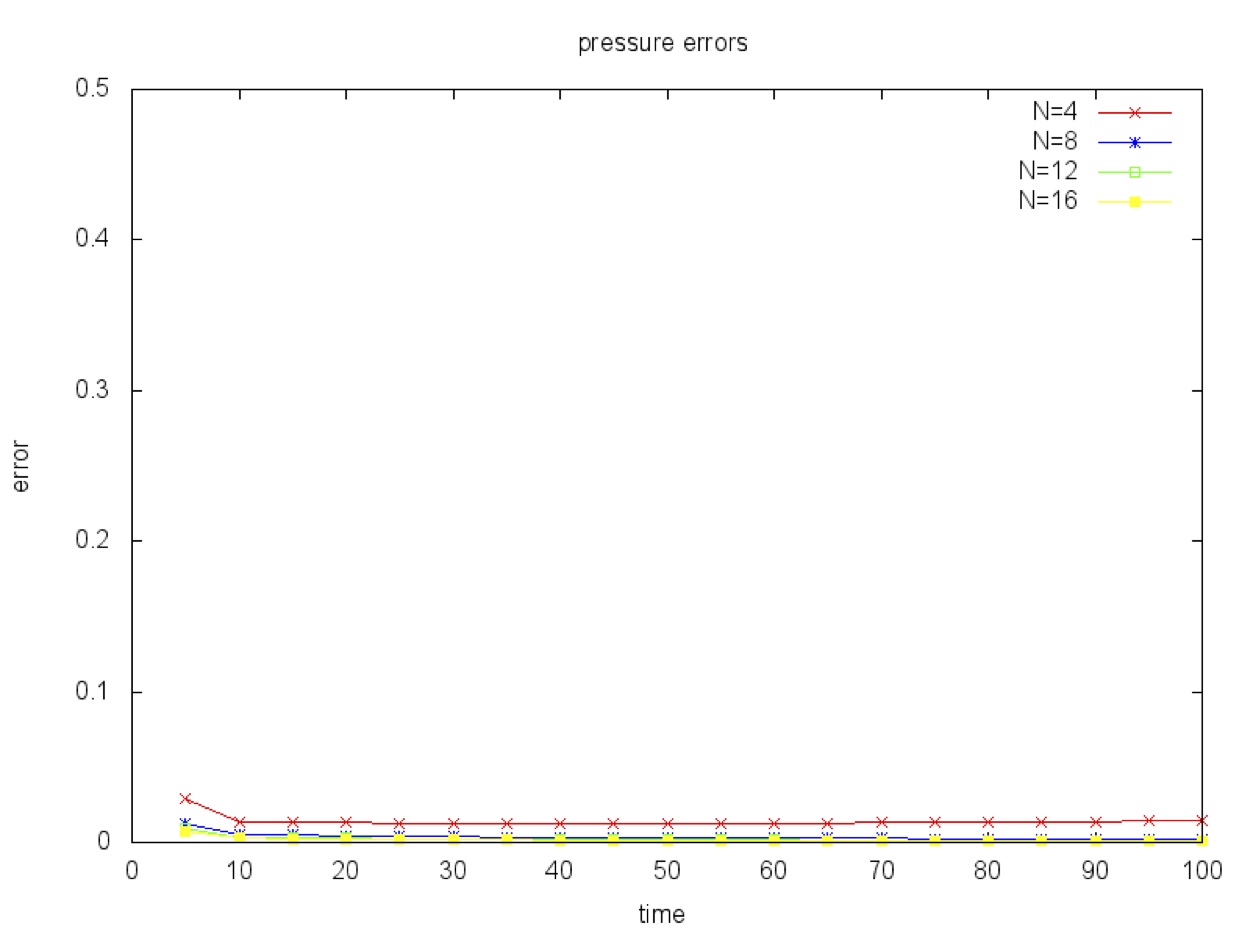} 
\includegraphics[width=0.45\linewidth]{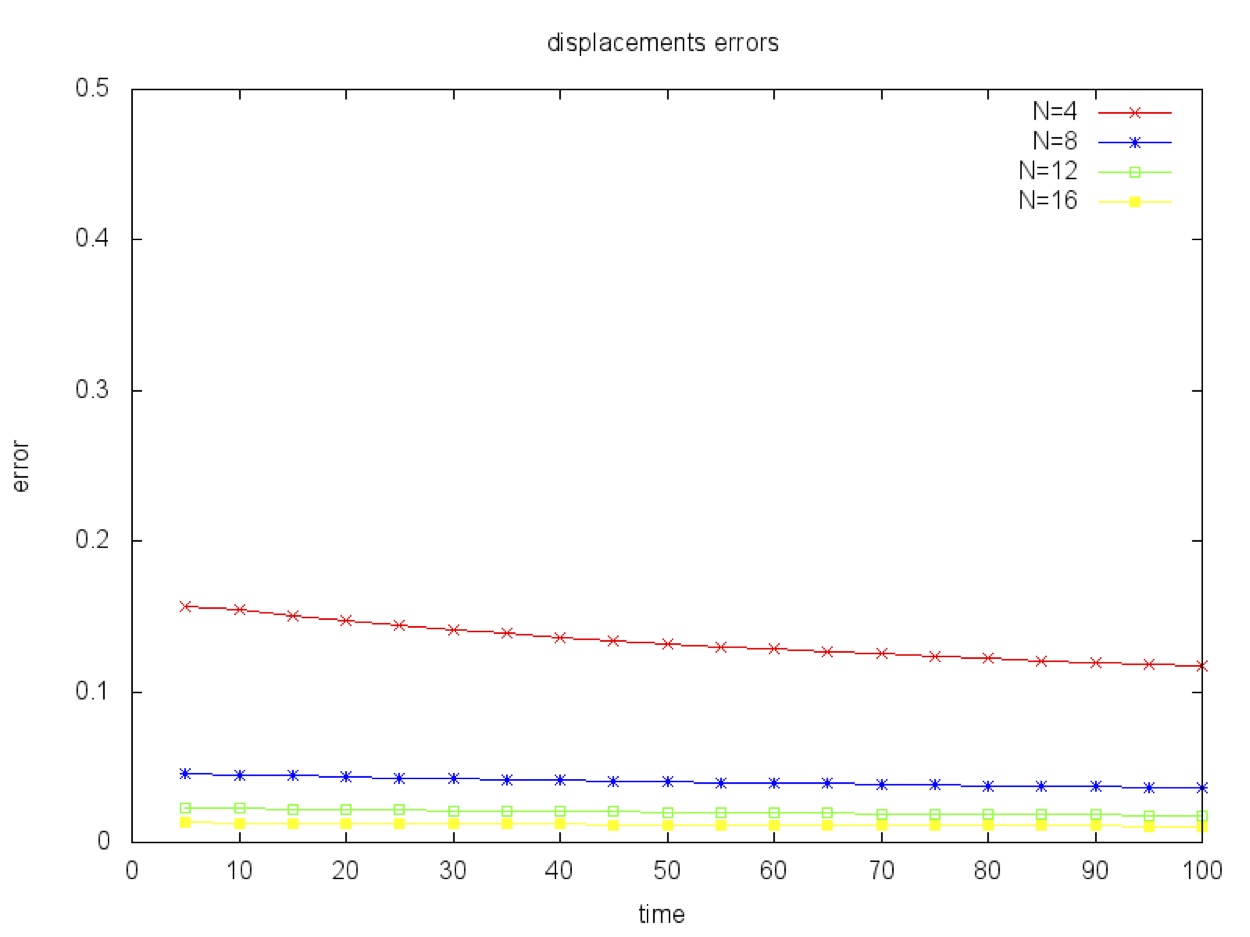} 
\end{center}
\caption{Weighted $L^2$ are on the top and $H^1$ are on the bottom. Errors for pressure are on the left and displacements are on the right  for Case 1 using the  fully coupled scheme.}
\label{fig:err-c1-c}
\end{figure}

\begin{figure}
\begin{center}
\includegraphics[width=0.45\linewidth]{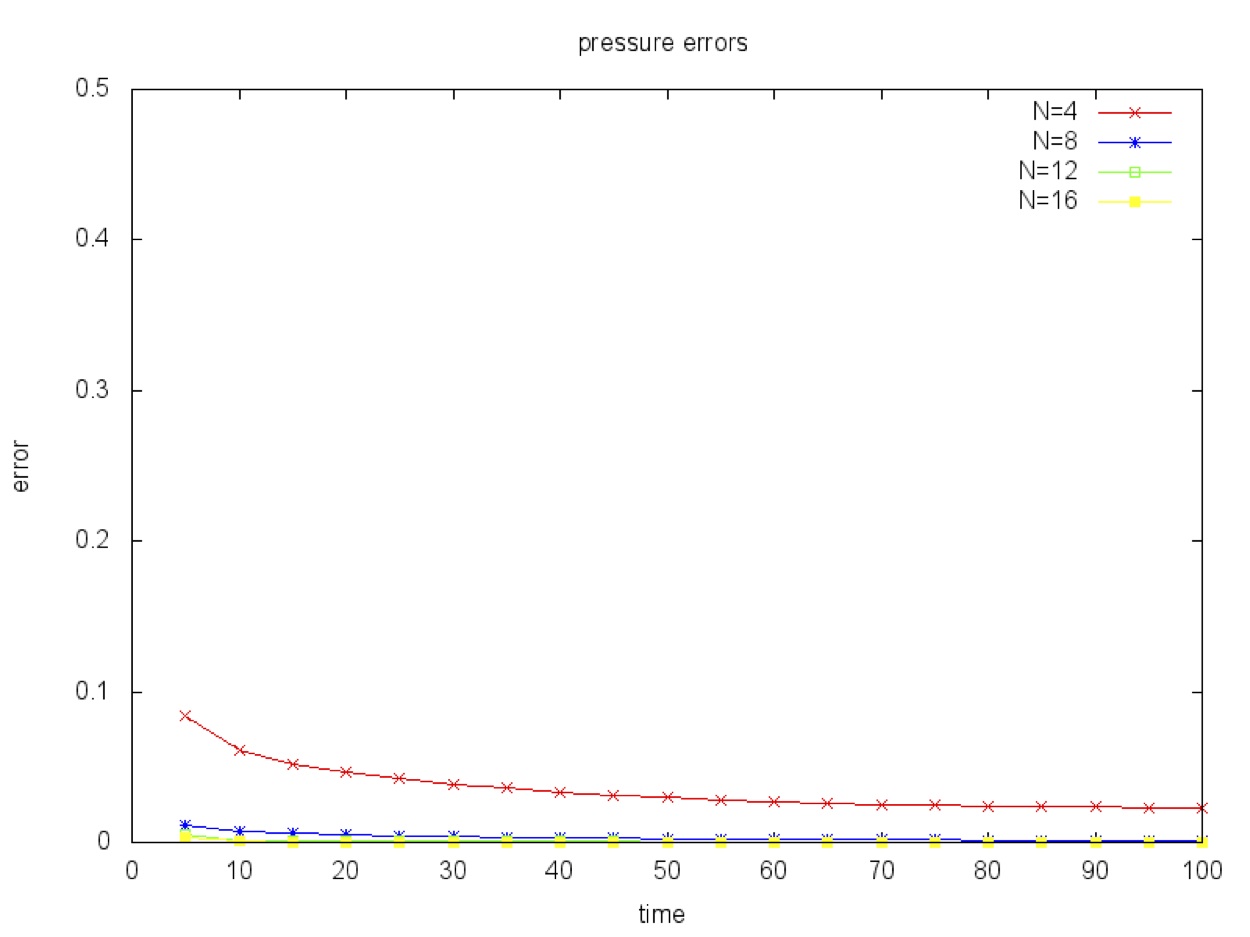} 
\includegraphics[width=0.45\linewidth]{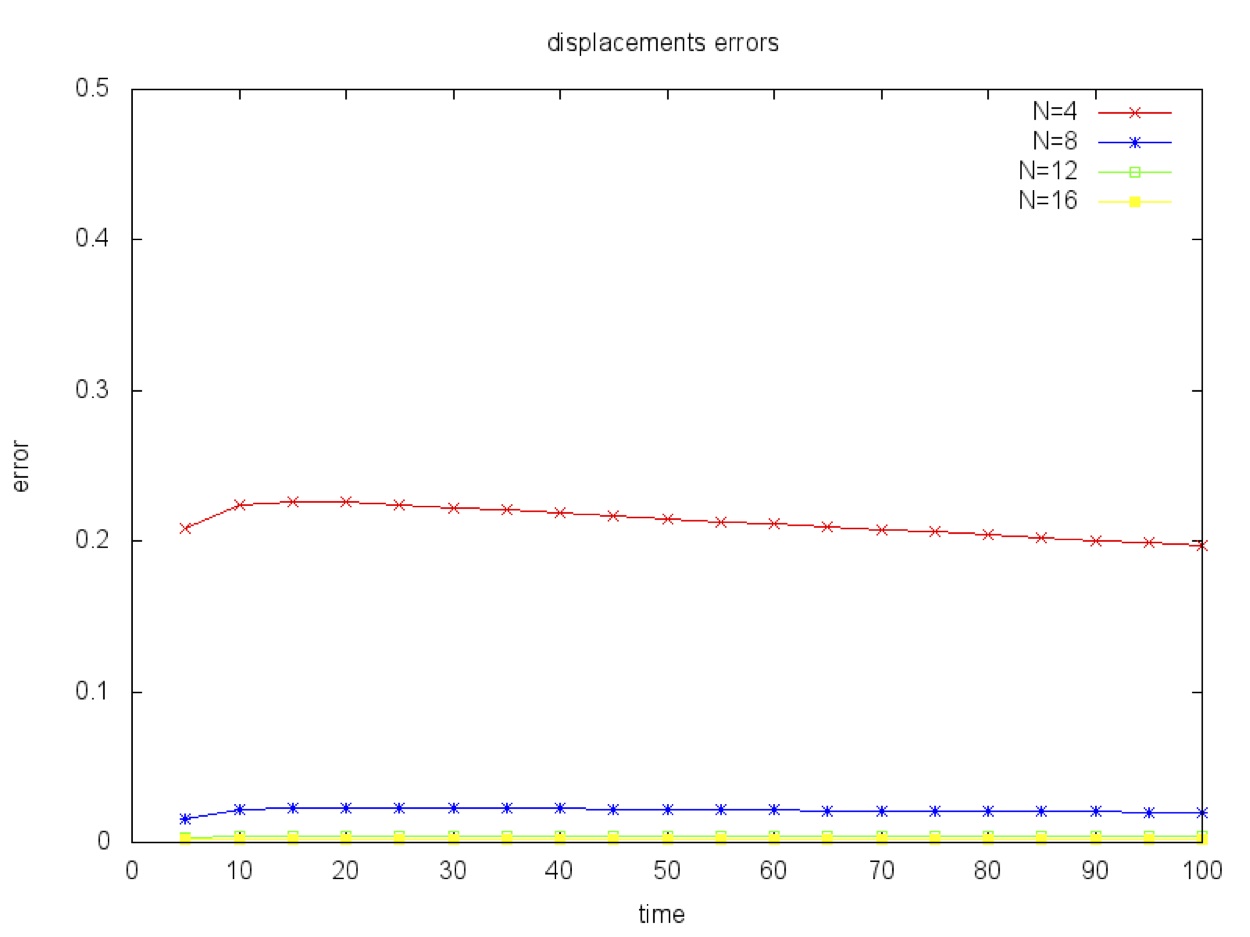} 
\\
\includegraphics[width=0.45\linewidth]{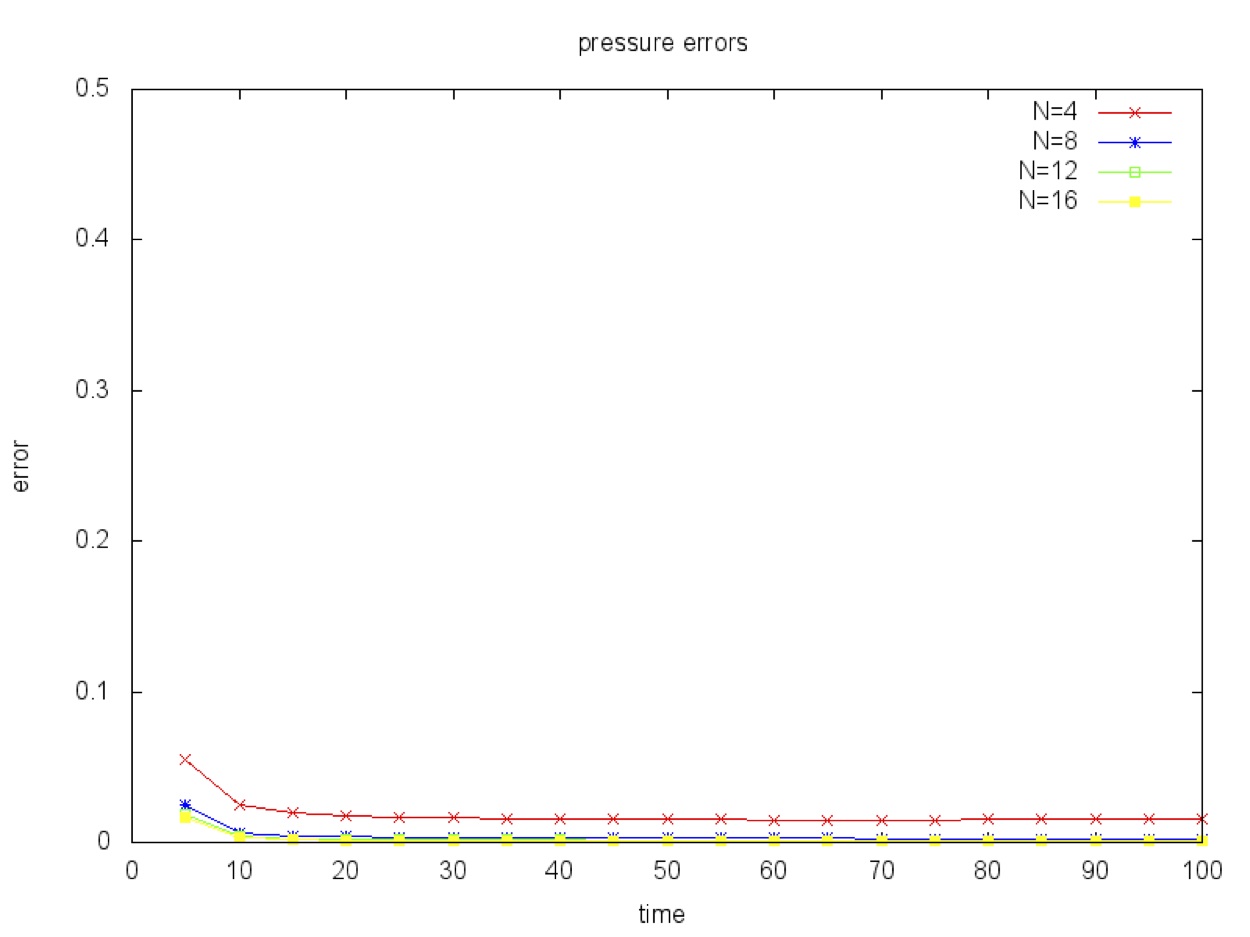} 
\includegraphics[width=0.45\linewidth]{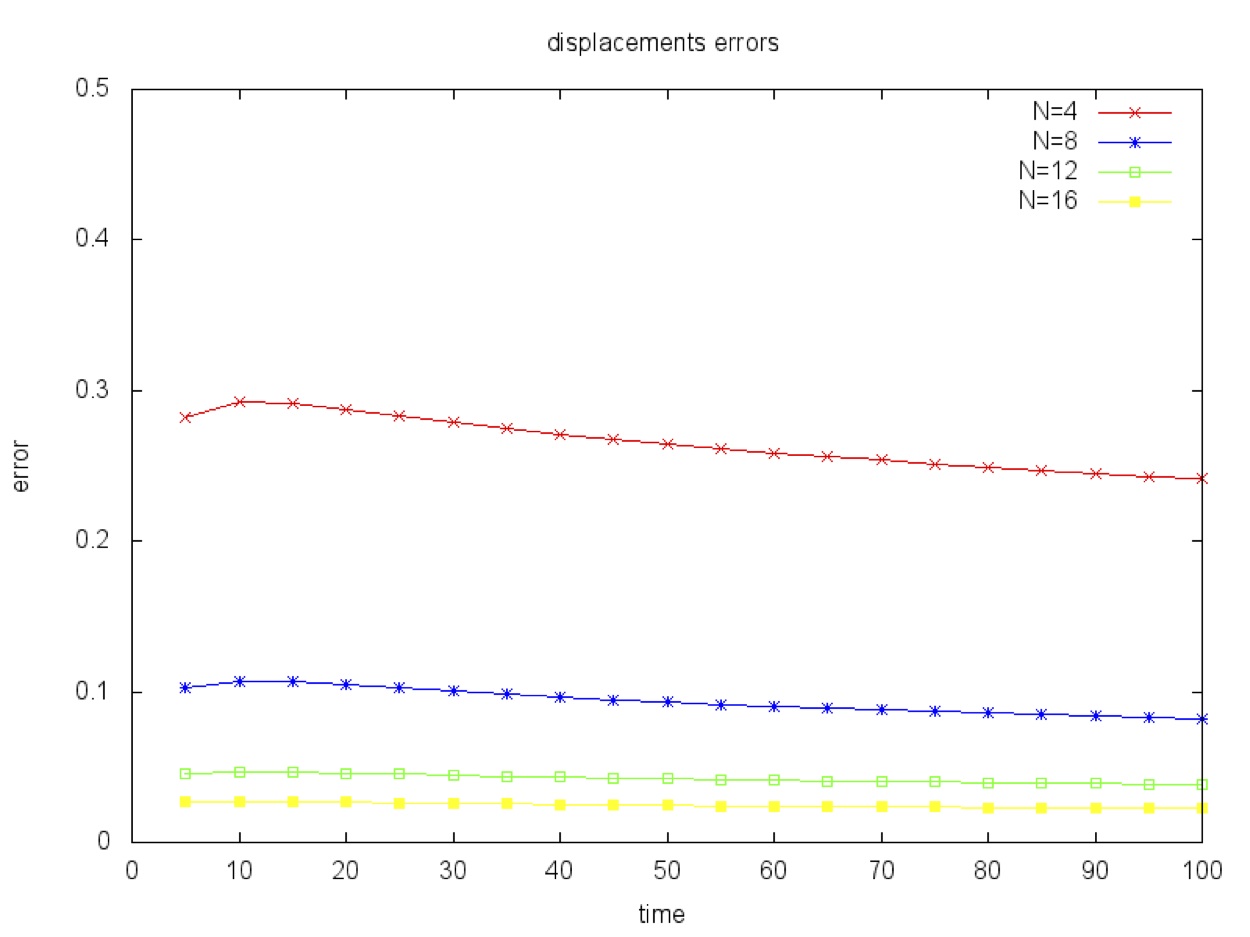} 
\end{center}
\caption{Weighted $L^2$ are on the top and $H^1$ are on the bottom. Errors for pressure are on the left and displacements are on the right  for Case 2 using the  fully coupled scheme.}
\label{fig:err-c2-c}
\end{figure}

We also include plots over time of the error with respect to number of basis functions used. We present the results from the fully coupled scheme. 
In Fig. \ref{fig:err-c1-c} and \ref{fig:err-c2-c} we show errors over time for $N_{\text{off}} = N_{\text{off}}^{p} = N_{\text{off}}^{u} = 4, 8, 12,$ and $16$ multiscale basis functions
 for each $\omega_{i}$
 Thus, the  dimensions of offline spaces are 432, 864, 1296 and 1728, respectively. 
 We observe that errors decrease as we increase the dimension of the offline space as expected.
We observe the errors in  Fig. \ref{fig:err-c1-c} are generally better than the errors Fig. \ref{fig:err-c2-c}, again, due to the lower contrast in Case 1. We see that in both cases most of the error vanished after the use of just 8 multiscale basis functions. In general, the error remains stable in time with a slight decrease over time.

\begin{table}[htp]
\begin{center}
\begin{tabular}[hp]{|c|c|cc|cc|}
\hline
 & & \multicolumn{2}{|c|}{Pressure errors}
 & \multicolumn{2}{|c|}{Displacements errors} \\
$N_{\text{off}}^p$  & dim($Q_{\text{off}},V_{\text{off}}$) & $L^2$ & $H^1$ & $L^2$ & $H^1$\\
\hline \hline
\multicolumn{6}{|c|}{ $N_{\text{off}}^u = 2$}  \\
\hline
2     & 216 & 0.06	&0.08   & 0.06	& 0.13 \\
\hline 
\multicolumn{6}{|c|}{$N_{\text{off}}^u = 4$}  \\
\hline
2     & 360 &  0.06	& 0.08   & 0.06	& 0.12 \\
4     & 432 &  0.01	& 0.01   & 0.04	& 0.11 \\
\hline 
\multicolumn{6}{|c|}{$N_{\text{off}}^u = 8$}  \\
\hline
2     & 648 &  0.06	& 0.08   & 0.02	& 0.06 \\
4     & 720 &  0.01	& 0.01   & 0.01	& 0.03 \\
8     & 864 &  0.0003	& 0.002  & 0.002	& 0.03 \\
\hline 
\multicolumn{6}{|c|}{$N_{\text{off}}^u= 12$}  \\
\hline
2     & 936 &  0.06	& 0.08   &  0.02	& 0.05 \\
4     & 1008 &  0.01	& 0.01   &  0.01	& 0.02 \\
8     & 1152 &  0.0003	& 0.002  &  0.0009	& 0.01 \\
12   & 1296 &  0.0001	& 0.001  &  0.0009	& 0.01 \\
\hline 
\multicolumn{6}{|c|}{$N_{\text{off}}^u = 16$}  \\
\hline
2     & 1224 &  0.06	& 0.08   &  0.02	& 0.05 \\
4     & 1296 &  0.01	& 0.01   &  0.01	& 0.01 \\
8     & 1440 &  0.0003	& 0.002  &  0.0008	& 0.01 \\
12   & 1584 &  0.0001	& 0.001  &  0.0007	& 0.01 \\
16   & 1728 &  0.0001	& 0.0007 &  0.0007 & 0.01 \\
\hline
\end{tabular}
\end{center}
\caption{Numerical results for Case 1 using the  fully coupled scheme.}
\label{tab:c1-c}
\end{table}

\begin{table}[htp]
\begin{center}
\begin{tabular}[hp]{|c|c|cc|cc|}
\hline
 & & \multicolumn{2}{|c|}{Pressure errors}
 & \multicolumn{2}{|c|}{Displacements errors} \\
$N_{\text{off}}^p$  & dim($Q_{\text{off}},V_{\text{off}}$) & $L^2$ & $H^1$ & $L^2$ & $H^1$\\
\hline \hline
\multicolumn{6}{|c|}{$N_{\text{off}}^u  = 2$}  \\
\hline
2     & 216 & 0.06	& 0.08  &  0.25	& 0.26 \\
\hline 
\multicolumn{6}{|c|}{$N_{\text{off}}^u = 4$}  \\
\hline
2     & 360 &  0.06	& 0.08  &  0.22	& 0.24 \\
4     & 432 &  0.02	& 0.01  &  0.19	& 0.24 \\
\hline 
\multicolumn{6}{|c|}{$N_{\text{off}}^u = 8$}  \\
\hline
2     & 648 &  0.06	& 0.08    &  0.08	& 0.13 \\
4     & 720 &  0.02	& 0.01    &  0.01	& 0.08 \\
8     & 864 &  0.001	& 0.002  &  0.01	& 0.08  \\
\hline 
\multicolumn{6}{|c|}{$N_{\text{off}}^u = 12$}  \\
\hline
2     & 936   &  0.06	    & 0.08    &  0.07   	& 0.11 \\
4     & 1008 &  0.02	    & 0.01    &  0.02	    & 0.04 \\
8     & 1152 &  0.0003	& 0.002  &  0.004	& 0.03 \\
12   & 1296 &  0.0001	& 0.001  &  0.004	& 0.03 \\
\hline 
\multicolumn{6}{|c|}{$N_{\text{off}}^u = 16$}  \\
\hline
2     & 1224 &  0.06	    & 0.08      &  0.07	& 0.11 \\
4     & 1296 &  0.02	    & 0.01      &  0.02	& 0.03 \\
8     & 1440 &  0.0003	& 0.002    &  0.001	& 0.02 \\
12   & 1584 &  0.0001	& 0.001    &  0.001	& 0.02 \\
16   & 1728 &  0.0001	& 0.0006  &  0.001	& 0.02 \\
\hline
\end{tabular}
\end{center}
\caption{Numerical results for Case 2 using the  fully coupled scheme. }
\label{tab:c2-c}
\end{table}

\begin{table}[htp]
\begin{center}
\begin{tabular}[hp]{|c|c|cc|cc|}
\hline
 & & \multicolumn{2}{|c|}{Pressure errors}
 & \multicolumn{2}{|c|}{Displacements errors} \\
$N_{\text{off}}^p$  & dim($Q_{\text{off}},V_{\text{off}}$) & $L^2$ & $H^1$ & $L^2$ & $H^1$\\
\hline \hline
\multicolumn{6}{|c|}{ $N_{\text{off}}^u = 2$}  \\
\hline
2     & 216 & 0.06	&0.08   & 0.06	& 0.13 \\
\hline 
\multicolumn{6}{|c|}{$N_{\text{off}}^u = 4$}  \\
\hline
2     & 360 &  0.06	& 0.08   & 0.06	& 0.12 \\
4     & 432 &  0.01	& 0.01   & 0.04	& 0.11 \\
\hline 
\multicolumn{6}{|c|}{$N_{\text{off}}^u = 8$}  \\
\hline
2     & 648 &  0.06	& 0.08   & 0.02	& 0.06 \\
4     & 720 &  0.01	& 0.01   & 0.01	& 0.03 \\
8     & 864 &  0.0003	& 0.002  & 0.002	& 0.03 \\
\hline 
\multicolumn{6}{|c|}{$N_{\text{off}}^u= 12$}  \\
\hline
2     & 936 &  0.06	& 0.08   &  0.02	& 0.05 \\
4     & 1008 &  0.01	& 0.01   &  0.01	& 0.02 \\
8     & 1152 &  0.0003	& 0.002  &  0.0009	& 0.01 \\
12   & 1296 &  0.0001	& 0.001  &  0.0009	& 0.01 \\
\hline 
\multicolumn{6}{|c|}{$N_{\text{off}}^u = 16$}  \\
\hline
2     & 1224 &  0.06	& 0.08   &  0.02	& 0.05 \\
4     & 1296 &  0.01	& 0.01   &  0.01	& 0.01 \\
8     & 1440 &  0.0003	& 0.002  &  0.0008	& 0.01 \\
12   & 1584 &  0.0001	& 0.001  &  0.0007	& 0.01 \\
16   & 1728 &  0.0001	& 0.0007 &  0.0007 & 0.01 \\
\hline
\end{tabular}
\end{center}
\caption{Numerical results  for Case 1 using the fixed-stress scheme.}
\label{tab:c1-s}
\end{table}

\begin{table}[htp]
\begin{center}
\begin{tabular}[hp]{|c|c|cc|cc|}
\hline
 & & \multicolumn{2}{|c|}{Pressure errors}
 & \multicolumn{2}{|c|}{Displacements errors} \\
$N_{\text{off}}^p$  & dim($Q_{\text{off}},V_{\text{off}}$)& $L^2$ & $H^1$ & $L^2$ & $H^1$\\
\hline \hline
\multicolumn{6}{|c|}{$N_{\text{off}}^u= 2$}  \\
\hline
2     & 216 & 0.30	& 0.26  &  0.45	& 0.46 \\
\hline 
\multicolumn{6}{|c|}{$N_{\text{off}}^u = 4$}  \\
\hline
2     & 360 &  0.30	& 0.26  &  0.42	& 0.45 \\
4     & 432 &  0.01	& 0.01  &  0.33	& 0.38 \\
\hline 
\multicolumn{6}{|c|}{$N_{\text{off}}^u = 8$}  \\
\hline
2     & 648 &  0.30	& 0.25    &  0.36	& 0.48 \\
4     & 720 &  0.006	& 0.01    &  0.04	& 0.15 \\
8     & 864 &  0.001	& 0.006  &  0.04	& 0.15  \\
\hline 
\multicolumn{6}{|c|}{$N_{\text{off}}^u = 12$}  \\
\hline
2     & 936   &  0.30	    & 0.25    &  0.37   	& 0.50 \\
4     & 1008 &  0.006   & 0.01    &  0.007	& 0.06 \\
8     & 1152 &  0.002	& 0.006  &  0.007	& 0.06 \\
12   & 1296 &  0.001	& 0.004  &  0.007	& 0.06 \\
\hline 
\multicolumn{6}{|c|}{$N_{\text{off}}^u = 16$}  \\
\hline
2     & 1224 &  0.30	    & 0.25      &  0.38	& 0.50 \\
4     & 1296 &  0.006   & 0.01      &  0.003	& 0.03 \\
8     & 1440 &  0.002	& 0.006    &  0.002	& 0.02 \\
12   & 1584 &  0.001	& 0.004    &  0.002	& 0.02 \\
16   & 1728 &  0.0009	& 0.003    &  0.002	& 0.02 \\
\hline
\end{tabular}
\end{center}
\caption{Numerical results  for Case 2 using the fixed-stress scheme.}
\label{tab:c2-s}
\end{table}

\subsection{GMsFEM with Randomized Oversampling}

In this section we consider the oversampling randomized algorithm proposed in \cite{randomized2014}. In this algorithm, instead of solving harmonic extensions (\ref{harmonic_ex} and \ref{harmonic_ex2}) for each fine grid node on the boundary, we solve a small number of harmonic extension local problems with  random boundary conditions. More precisely, we let 
\[
\psi_{j}^{\omega_i, \text{snap}}=r_j, \quad x \in \partial\omega_i,
\]
where $r_j$ are independent identical distributed standard Gaussian random vectors on the fine grid nodes of the boundary. 
The advantage of this algorithm lies in the fact that a much fewer number of snapshot basis functions are calculated, while maintaining accuracy. 
In addition, we will use an oversampling strategy. This is done to reduce the mismatching effects of boundary conditions imposed artificially in the construction of snapshot basis functions. 
We will denote the extended coarse grid neighborhood for $t=1,2,\dots$, by $\omega^+_i=\omega_{i}+t$. Here for example, $\omega^+_i=\omega_{i}+1$, would mean the coarse grid neighborhood plus all 1 layer of adjacent fine grida of $\omega_{i}$, and so on.

\begin{figure}
\begin{center}
\includegraphics[width=0.45\linewidth]{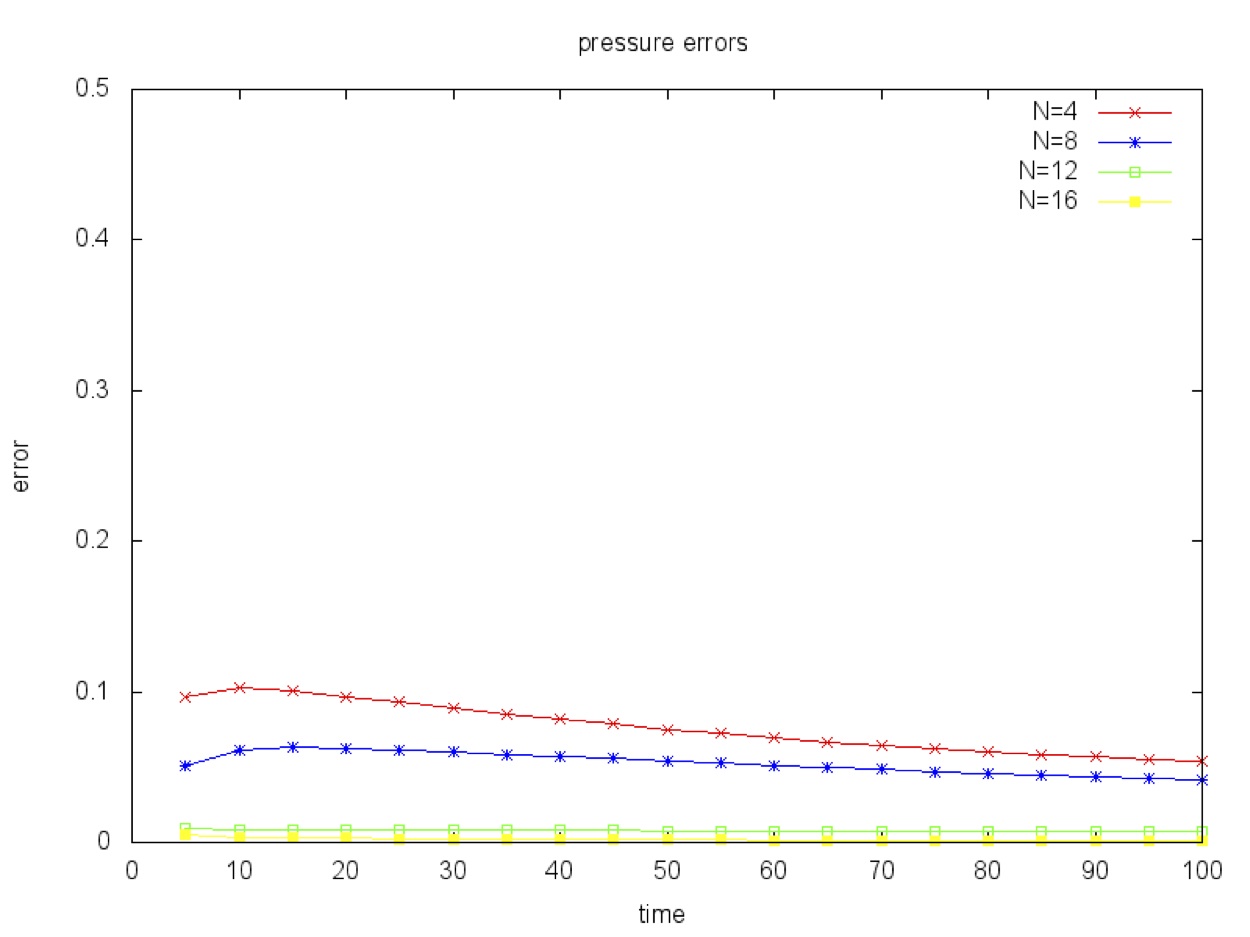} 
\includegraphics[width=0.45\linewidth]{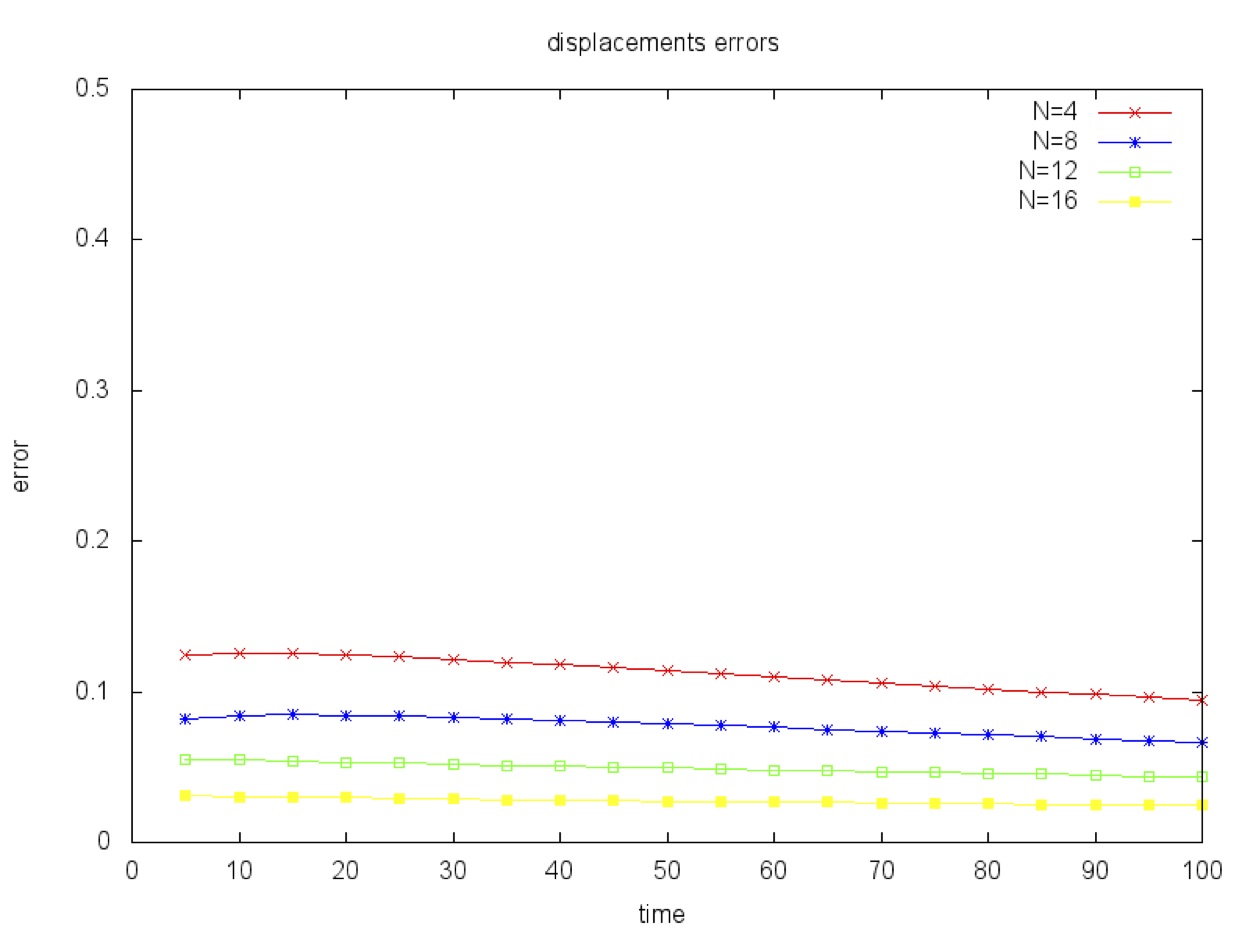} 
\\
\includegraphics[width=0.45\linewidth]{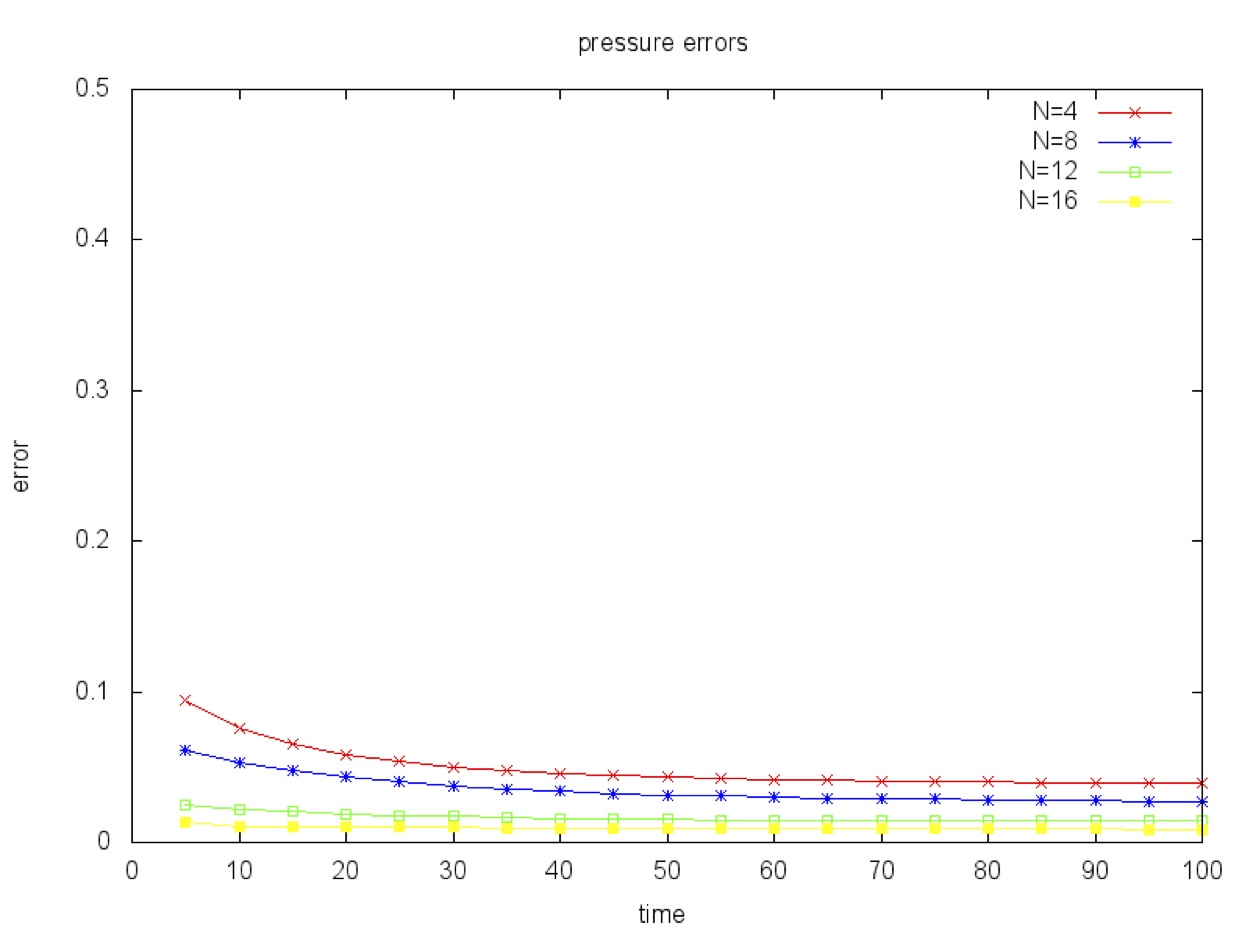} 
\includegraphics[width=0.45\linewidth]{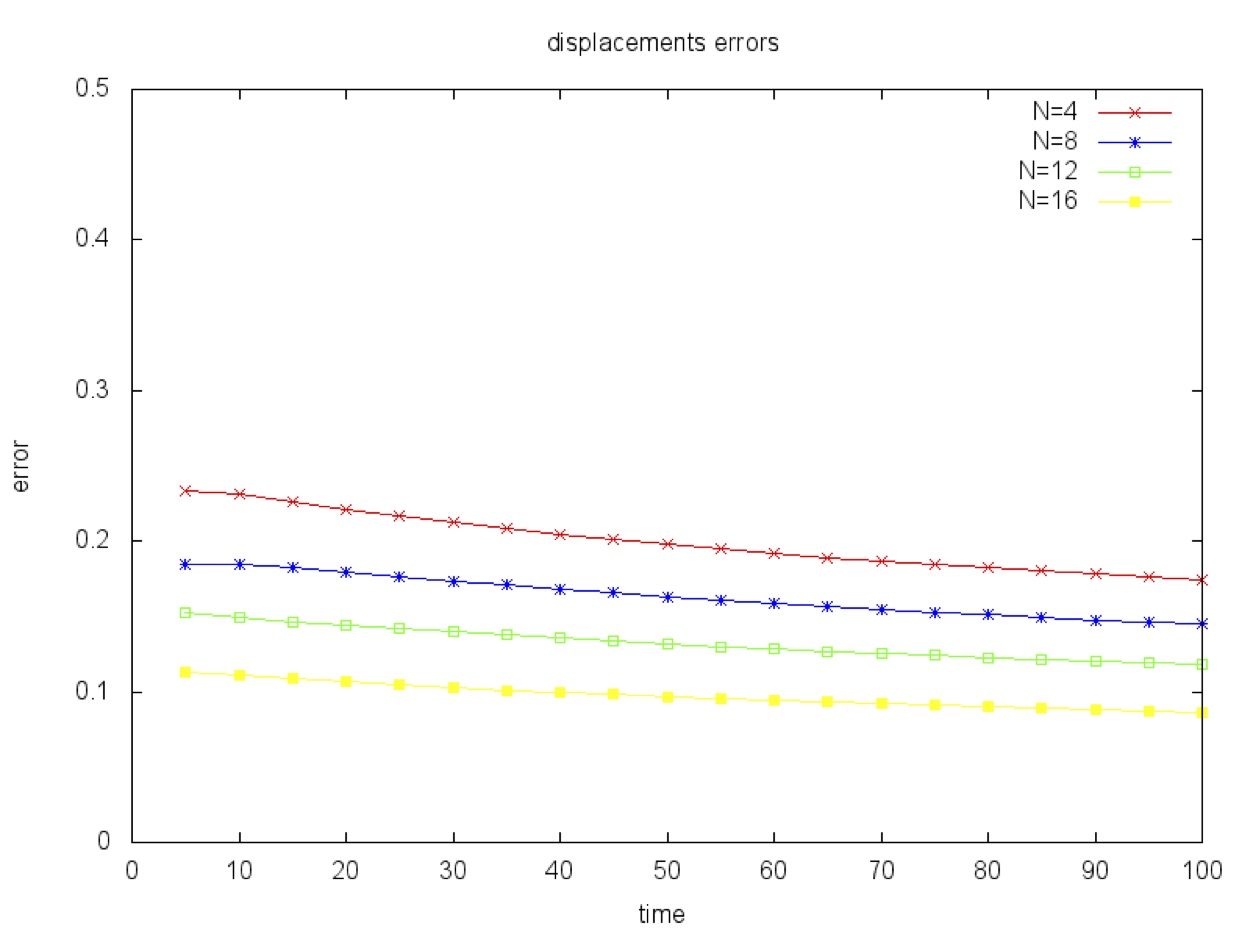} 
\end{center}
\caption{Weighted $L^2$ are on the top and $H^1$ are on the bottom. Errors for pressure are on the left and displacements are on the right  for Case 1 using randomized GMsFEM with oversampling, $\omega^+_i = \omega_i + 4$.}
\label{pic:err-c1-r}
\end{figure}

\begin{figure}
\begin{center}
\includegraphics[width=0.45\linewidth]{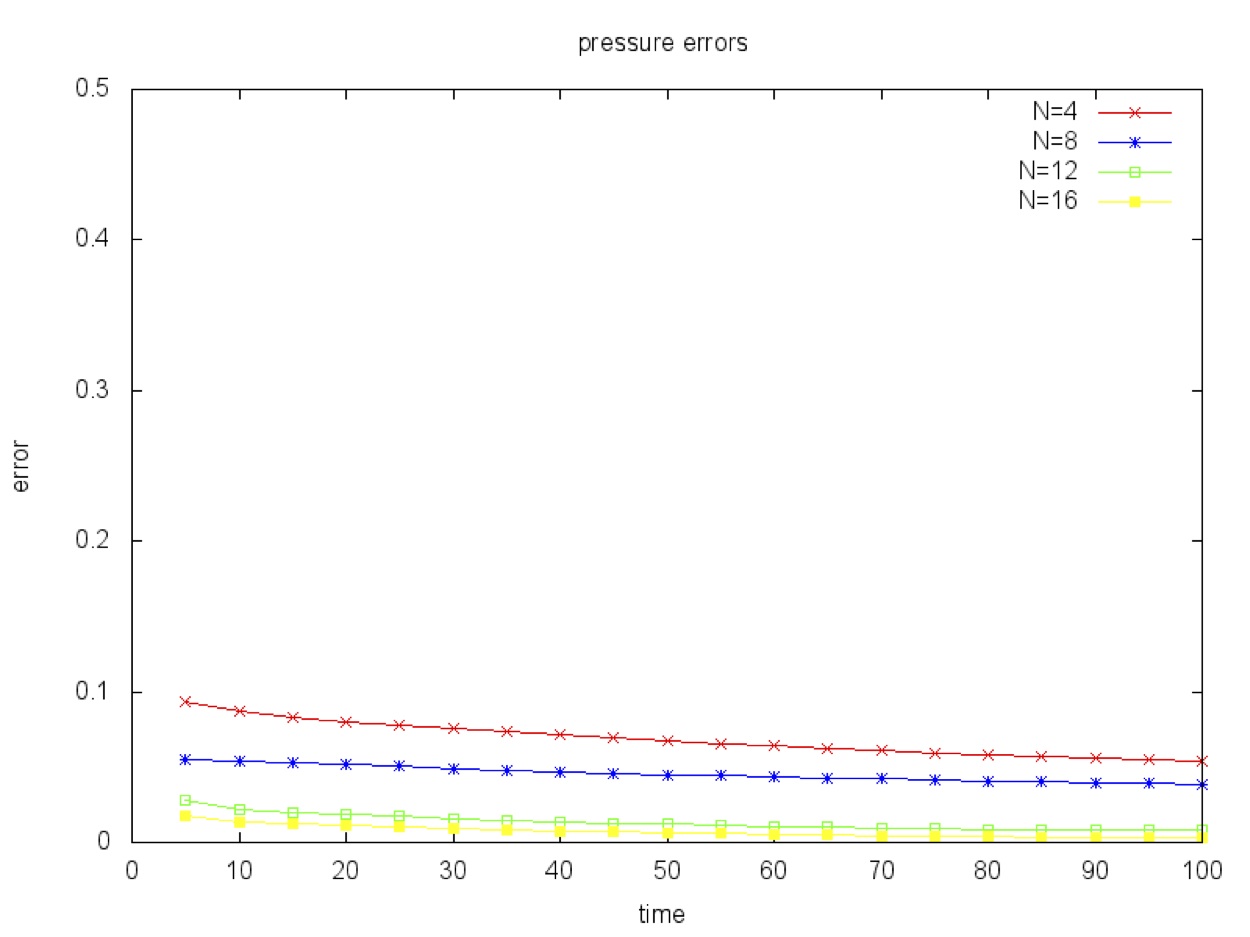} 
\includegraphics[width=0.45\linewidth]{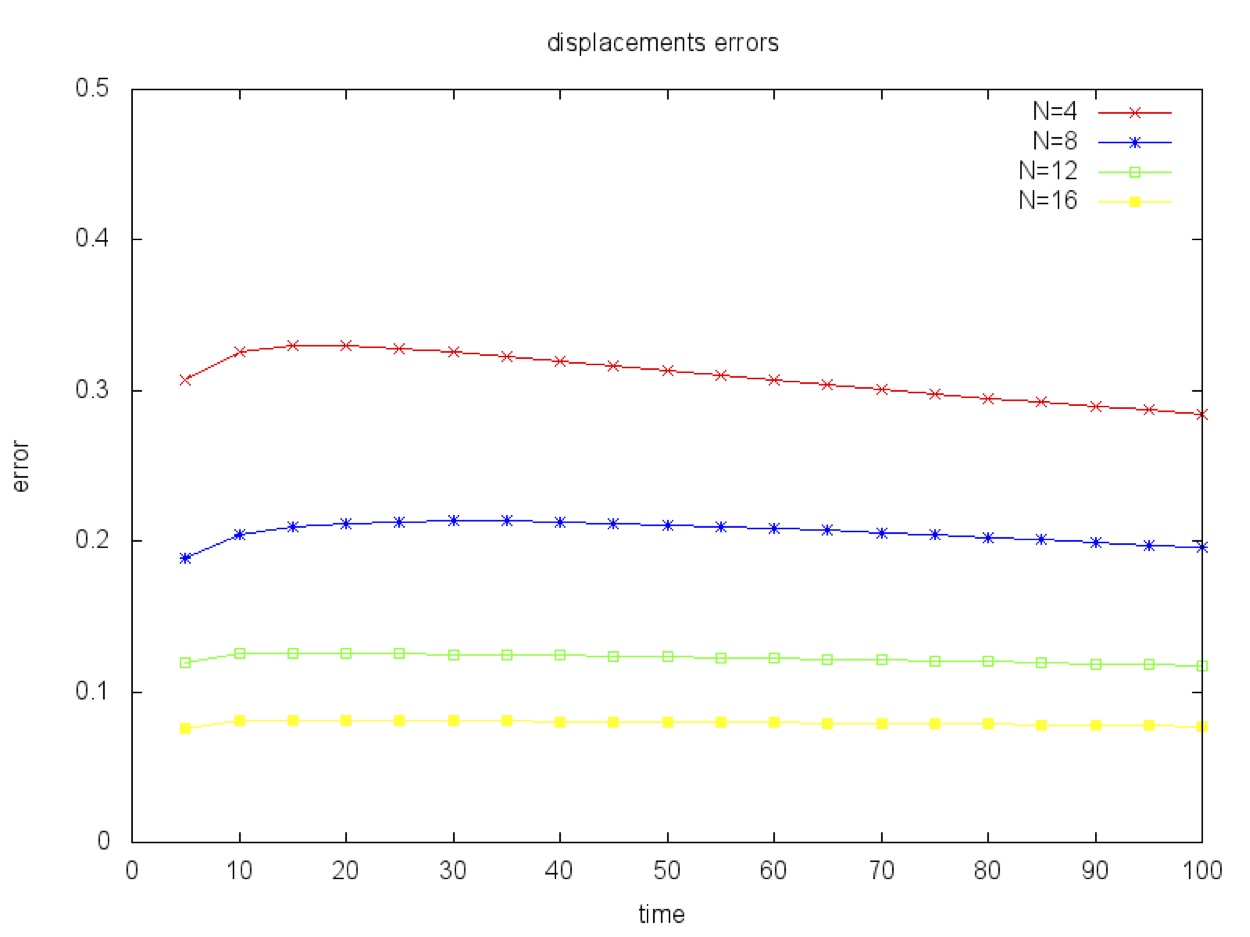} 
\\
\includegraphics[width=0.45\linewidth]{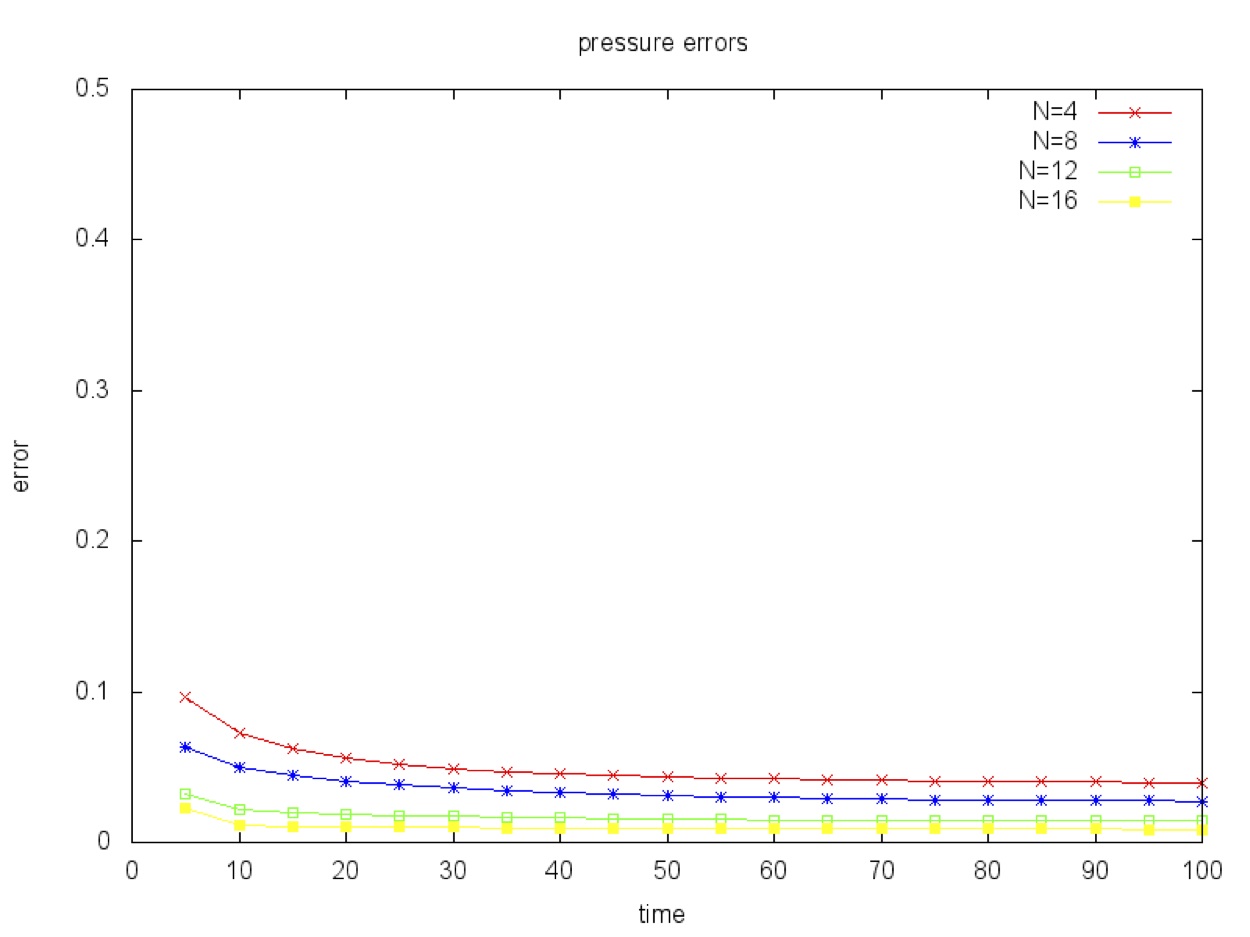} 
\includegraphics[width=0.45\linewidth]{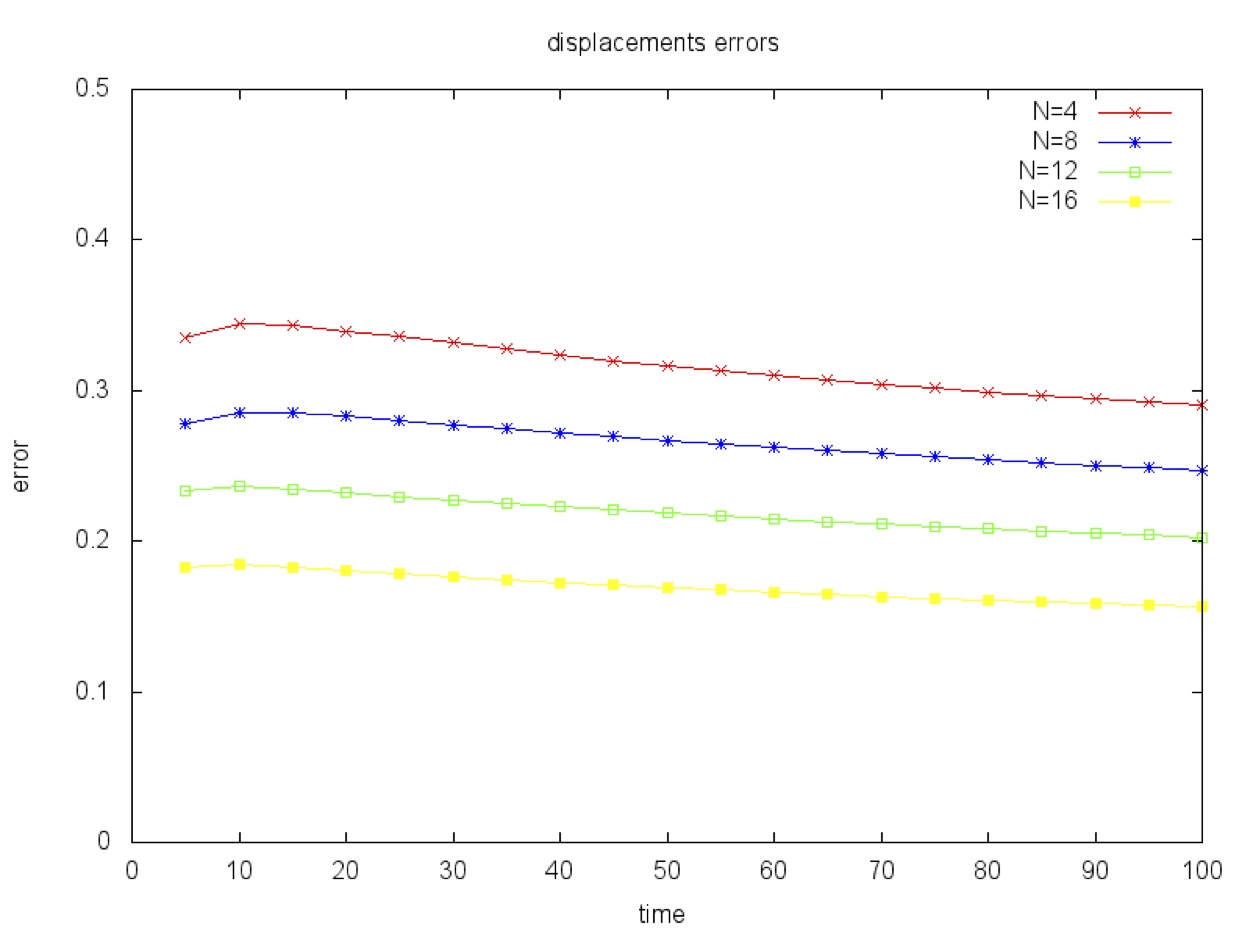} 
\end{center}
\caption{Weighted $L^2$ are on the top and $H^1$ are on the bottom. Errors for pressure are on the left and displacements are on the right  for Case 2  using randomized GMsFEM with oversampling, $\omega^+_i = \omega_i + 4$.}
\label{pic:err-c2-r}
\end{figure}

In Fig. \ref{pic:err-c1-r} and Fig. \ref{pic:err-c2-r}, we show the weighted $L^2$ and $H^1$ errors  over time  for Case 1 and 2 using the randomized GMsFEM with oversampling using different numbers of multiscale basis functions.
The oversampled region $\omega^+_i = \omega_i + 4$  is chosen, that is,  the oversampled region contains an extra 4 fine grid cells layers around $\omega_i$. 
Here, we use only the fully coupled scheme.  
We use a snapshot ratio of $36\%$ between the standard number of snapshots and the randomized algorithm. 
Comparing results from Fig. \ref{pic:err-c1-r}, the randomized algorithm,  to Fig. \ref{fig:err-c1-c}, the standard GMsFEM,  we observe that the randomized algorithm is slightly less accurate but at the 
advantage of having less snapshot solutions required. 

In Table \ref{tab:c1-r} and \ref{tab:c2-r} we investigate  the effect of the oversampling $\omega^+_i = \omega_i + t$ as we increase the number of fine grid extensions for $t = 0, 2, 4$ and $6$. 
 We present the data of the randomized snapshots for last time step. We see that 
 oversampling can help to improve the results initially, 
 but the improvements level off as large oversampling domains do not give significant improvement in the solution accuracy. Again the effects of the high contrast of Case 2 can be seen in the data
 as the oversampling performs slightly worse than in the lower contrast regime.

\begin{table}[htp]
\begin{center}
\begin{tabular}[hp]{|c|cc|cc|}
\hline
 & \multicolumn{2}{|c|}{pressure errors}
 & \multicolumn{2}{|c|}{displacements errors} \\
$N_{\text{off}}$  & $L^2$ & $H^1$ & $L^2$ & $H^1$\\
\hline \hline 
\multicolumn{5}{|c|}{without oversampling, $\omega^+_i = \omega_i$}  \\
\hline
4       &  0.05   & 0.04  &  0.16   	& 0.22 \\
8       &  0.05   & 0.03  &  0.16	    & 0.21 \\
12     &  0.02	& 0.02  &  0.16	& 0.21 \\
16     &  0.004	& 0.01  &  0.15	& 0.21 \\
\hline 
\multicolumn{5}{|c|}{with oversampling, $\omega^+_i = \omega_i + 2$}  \\
\hline
4       &  0.05   & 0.03    &  0.14  & 0.21 \\
8       &  0.04   & 0.03    &  0.12	& 0.19 \\
12     &  0.007	& 0.01    &  0.09	& 0.17 \\
16     &  0.002	& 0.009  &  0.08	& 0.16 \\
\hline 
\multicolumn{5}{|c|}{with oversampling, $\omega^+_i = \omega_i + 4$}  \\
\hline
4       &  0.05   & 0.03    &  0.09  & 0.17 \\
8       &  0.04   & 0.02    &  0.06	& 0.14 \\
12     &  0.006	& 0.01    &  0.04	& 0.11 \\
16     &  0.001	& 0.008  &  0.02	& 0.08 \\
\hline  
\multicolumn{5}{|c|}{with oversampling, $\omega^+_i = \omega_i + 6$}  \\
\hline
4       &  0.05   & 0.03    &  0.09  & 0.17 \\
8       &  0.04   & 0.02    &  0.06	& 0.13 \\
12     &  0.009	& 0.01    &  0.02	& 0.09 \\
16     &  0.002	& 0.007  &  0.02	& 0.07 \\
\hline 
\end{tabular}
\end{center}
\caption{Numerical tests for Case 1 using randomized GMsFEM with and without oversampling for $N_{\text{off}} =  N_{\text{off}}^u = N_{\text{off}}^p$}
\label{tab:c1-r}
\end{table}

\begin{table}[htp]
\begin{center}
\begin{tabular}[hp]{|c|cc|cc|}
\hline
 & \multicolumn{2}{|c|}{pressure errors}
 & \multicolumn{2}{|c|}{displacements errors} \\
$N$  & $L^2$ & $H^1$ & $L^2$ & $H^1$\\
\hline \hline 
\multicolumn{5}{|c|}{without oversampling, $\omega^+_i = \omega_i$}  \\
\hline
4       &  0.05   & 0.04  &  0.36   	& 0.31 \\
8       &  0.05   & 0.03  &  0.35	    & 0.31 \\
12     &  0.02	& 0.02  &  0.34	& 0.31 \\
16     &  0.006	& 0.01  &  0.34	& 0.31 \\
\hline 
\multicolumn{5}{|c|}{with oversampling, $\omega^+_i = \omega_i + 2$}  \\
\hline
4       &  0.05   & 0.03    &  0.33  & 0.31 \\
8       &  0.04   & 0.03    &  0.30	& 0.30 \\
12     &  0.01	& 0.01    &  0.25	& 0.27 \\
16     &  0.009	& 0.009  &  0.22	& 0.25 \\
\hline 
\multicolumn{5}{|c|}{with oversampling, $\omega^+_i = \omega_i + 4$}  \\
\hline
4       &  0.05   & 0.03    &  0.28  & 0.29 \\
8       &  0.03   & 0.02    &  0.19	& 0.24 \\
12     &  0.007	& 0.01    &  0.11	& 0.20 \\
16     &  0.002	& 0.008  &  0.07	& 0.15 \\
\hline  
\multicolumn{5}{|c|}{with oversampling, $\omega^+_i = \omega_i + 6$}  \\
\hline
4       &  0.05   & 0.03    &  0.24  & 0.27 \\
8       &  0.04   & 0.02    &  0.14	& 0.22 \\
12     &  0.01	& 0.01    &  0.07	& 0.17 \\
16     &  0.002	& 0.007  &  0.06	& 0.14 \\
\hline 
\end{tabular}
\end{center}
\caption{Numerical tests for Case 2 using randomized GMsFEM with and without oversampling for $N_{\text{off}} =  N_{\text{off}}^u = N_{\text{off}}^p$}
\label{tab:c2-r}
\end{table}

\section{Conclusion} 
Simulating poroelasticity is difficult due the complex heterogeneities and  because of the complexity of gridding the flow and mechanics regimes in such media. %
Therefore, in this paper we developed a Generalized Multiscale Finite Element Method for a linear poroelastic media. We first presented the general poroelasticity framework of Biot and its subsequent solution by fixed stress time splitting methods. Although fully coupled schemes are considered numerically, this splitting lays the framework for the application of the GMsFEM to the decoupled poroelastic equations. We then outline the construction of the multiscale spaces in both fluid and mechanics regimes. The algorithm is then implemented on a single geometry with two different cases of elastic parameters. We show the errors relative to the fine scale solution over time and with varying multiscale basis functions. Finally, we implemented  oversampling strategies and randomized boundary conditions when solving for the  snapshot space. As in cases of reservoir compaction, the permeability may depend on pressure resulting in a nonlinear relation. In future studies, we will develop a GMsFEM for such nonlinear poroelastic problems.

\end{document}